\documentclass[12pt]{article}
\usepackage{amsmath,amsthm,amssymb,amscd}
\theoremstyle{plain}
\newtheorem{thmi}{Theorem}
\newtheorem{theorem}{Theorem}[section]
\newtheorem{lemma}[theorem]{Lemma}
\newtheorem{proposition}[theorem]{Proposition}
\newtheorem{corollary}[theorem]{Corollary}
\theoremstyle{definition}
\newtheorem{definition}[theorem]{Definition}
\theoremstyle{remark}
\newtheorem{remark}[theorem]{Remark}
\numberwithin{equation}{section}
\DeclareMathOperator{\Spec}{Spec}
\DeclareMathOperator{\ord}{ord}
\DeclareMathOperator{\sgn}{sgn}
\DeclareMathOperator{\Res}{Res}
\DeclareMathOperator{\Hom}{Hom}
\DeclareMathOperator{\Ext}{Ext}
\DeclareMathOperator{\Cl}{Cl}

\DeclareMathOperator{\Ker}{Ker}
\DeclareMathOperator{\Coker}{Coker}

\title{On the Hodge structure of degenerating hypersurfaces in toric
varieties}
\author{Atsushi Ikeda}
\date{}
\begin{document}
\maketitle
\begin{abstract}
 We introduce an algebraic method for describing the Hodge filtration of
 degenerating hypersurfaces in projective toric varieties.
 For this purpose, we show some fundamental properties of logarithmic
 differential forms on proper equivariant morphisms of toric varieties.
\end{abstract}
\section{Introduction}
The method for describing Hodge filtration by using polynomials
originated in the study of hypersurfaces in projective spaces by
Griffiths \cite{g}.
The theory is extended to the case on hypersurfaces in simplicial
projective toric varieties by Dolgachev \cite{do}, Steenbrink \cite{st},
Batyrev and Cox \cite{bc}.
The purpose of this paper is to apply their idea to degenerating
families of hypersurfaces in projective toric varieties.\par
Let $\pi:\mathbf{P}\rightarrow\mathbf{A}$ be a proper surjective
equivariant morphism of toric varieties over an algebraically closed
field $k$.
Here we assume for simplicity that $\mathbf{P}$ is nonsingular,
$\mathbf{A}$ is an affine space
$\mathbf{A}=\Spec{k[t_{1},\dots,t_{m}]}$, and the characteristic of $k$
is $0$.
Let $X$ be a hypersurface in $\mathbf{P}$.
In the case when $\pi$ is flat and geometrically connected, $\pi$ gives
a trivial fibration of a nonsingular complete toric variety over the
open torus of $\mathbf{A}$, and degenerated fibers appear at the outside
of the open torus.
Hence we can consider $X\rightarrow\mathbf{A}$ to be a degenerating
family of hypersurfaces in the complete toric variety.
We define a Jacobian ring for the family, and describe the Hodge
filtration of the family by using the Jacobian ring.\par
We denote by $D_{1},\dots, D_{s}$ all prime divisors invariant under the
torus action on $\mathbf{P}$.
The homogeneous coordinate ring of $\mathbf{P}$ is defined in \cite{c}
as a polynomial ring $S_{\mathbf{P}}=k[z_{1},\cdots,z_{s}]$ which has a
grading valued in the divisor class group $\Cl{(\mathbf{P})}$;
$$
\deg{z_{i}}=[D_{i}]\in\Cl{(\mathbf{P})}.
$$
We can assume that $\pi(D_{i})=\mathbf{A}$ for $1\leq i\leq r$, and
$\pi(D_{j})$ is contained in the divisor $\{t_{1}\cdots t_{m}=0\}$ for
$r+1\leq j\leq s$.
The hypersurface $X$ is defined by a $\Cl{(\mathbf{P})}$-homogeneous
polynomial $F\in S_{\mathbf{P}}$.
Then we define the Jacobian ring of $X$ over $\mathbf{A}$ by
$$
R_{X/\mathbf{A}}=
S_{\mathbf{P}}/
(\frac{\partial F}{\partial z_{1}},\dots,
\frac{\partial F}{\partial z_{r}}),
$$
which is a $\Cl{(\mathbf{P})}$-graded $k[t_{1},\dots,t_{m}]$-algebra.
For $\beta\in\Cl{(\mathbf{P})}$, the degree $\beta$ part of
$R_{X/\mathbf{A}}$ is denoted by $R_{X/\mathbf{A}}^{\beta}$, which is a
finitely generated $k[t_{1},\dots,t_{m}]$-module.\par
The Hodge filtration of the degenerating family is defined by using the
sheaf of relative logarithmic differential forms, so we consider the
situation above with logarithmic structure \cite{k}.
We define a logarithmic structure on $\mathbf{P}$ by the divisor
$E=\sum_{j=r+1}^{s}D_{j}$, and define a logarithmic structure on
$\mathbf{A}$ by the divisor $\{t_{1}\cdots t_{m}=0\}$.
Then $\pi$ is a logarithmically smooth morphism, and the logarithmic
structure of the general fiber is trivial.
The sheaf of relative logarithmic differential $p$-forms is denoted by
$\omega_{\mathbf{P}/\mathbf{A}}^{p}$, which is a locally free
$\mathcal{O}_{\mathbf{P}}$-modules.
We define a logarithmic structure on $X$ by the restriction of the
logarithmic structure on $\mathbf{P}$.
The next theorem is our main result, where we need not assume that $\pi$
is flat and geometrically connected.
\begin{thmi}[Theorem~$\ref{main}$]\label{imain}
 If $X$ is ample and logarithmic smooth over an affine open subvariety
 $U=\Spec{A_{U}}$ of $\mathbf{A}$, then for $0\leq p\leq n-1$, there is
 a natural isomorphism of $A_{U}$-modules
 $$
 H^{n-p-1}(\mathbf{P}_{U},
 \omega_{\mathbf{P}/\mathbf{A}}^{p+1}(\log{X}))
 \simeq
 A_{U}\otimes_{k[t_{1},\dots,t_{m}]}R_{X/\mathbf{A}}^{[(n-p)X-D]},
 $$
 where $n=\dim{\mathbf{P}}-\dim{\mathbf{A}}$,
 $\mathbf{P}_{U}=U\times_{\mathbf{A}}\mathbf{P}$ and
 $D=\sum_{i=1}^{r}D_{i}$.
\end{thmi}
If $m=0$, then $\mathbf{P}$ is a nonsingular complete variety, and the
logarithmic structure on $\mathbf{P}$ is trivial.
For an ample smooth hypersurface $X$ in $\mathbf{P}$, the isomorphism in
Theorem~$\ref{imain}$ is
$$
H^{n-p-1}(\mathbf{P},
\varOmega_{\mathbf{P}}^{p+1}(\log{X}))
\simeq
R_{X}^{[(n-p)X-D]}.
$$
This is proved in \cite{bc}, so Theorem~$\ref{imain}$ is a
generalization of the result by Batyrev and Cox.\par
In the case when
$\mathbf{P}\rightarrow\mathbf{A}^{1}\times\mathbf{P}^{n}$ is the blowing
up along a point, and $\pi$ is defined by the composition with the first
projection to $\mathbf{A}^{1}$, the logarithmically smooth family
$X_{U}\rightarrow U$ is a semistable degeneration of hypersurfaces in
$\mathbf{P}^{n}$.
This example is studied by Saito in \cite{sa}, which is the first work
for describing the Hodge filtration of degenerating hypersurfaces by
using the Jacobian rings.\par
The key of the proof of Theorem~$\ref{imain}$ is the following two
fundamental property for the sheaf of relative logarithmic differential
forms on $\mathbf{P}$.
The first property is a generalization of the Bott vanishing theorem:
\begin{thmi}[Corollary~$\ref{v}$]\label{iv}
 If $\mathcal{L}$ is an ample invertible sheaf on $\mathbf{P}$, then for
 $p\geq0$ and $q\geq1$,
 $$
 H^{q}(\mathbf{P},
 \omega_{\mathbf{P}/\mathbf{A}}^{p}
 \otimes_{\mathcal{O}_{\mathbf{P}}}\mathcal{L})=0.
 $$
\end{thmi}
The second property is a generalization of the Euler exact sequence:
\begin{thmi}[Theorem~$\ref{eu}$]\label{ie}
 There is an exact sequence of $\mathcal{O}_{\mathbf{P}}$-modules
 $$
 0\longrightarrow
 \omega_{\mathbf{P}/\mathbf{A}}^{1}
 \longrightarrow
 \bigoplus_{i=1}^{r}
 \mathcal{O}_{\mathbf{P}}(-D_{i})
 \longrightarrow
 \mathcal{O}_{\mathbf{P}}\otimes_{\mathbf{Z}}
 \Cl{(\mathbf{P}\smallsetminus E)}
 \longrightarrow0.
 $$
\end{thmi}
We prove Theorem~$\ref{iv}$ and Theorem~$\ref{ie}$ by using the
Poincar\'{e} residue map for the sheaf of relative logarithmic
differential forms, that is the idea of Batyrev and Cox \cite{bc}.\par
This paper proceed as follows.
In Section~$\ref{s2}$, we consider invertible sheaves on a toric variety
with a proper equivariant morphism to an affine toric variety, and we
give a characterization of base point freeness or ampleness of the
linear system by the terminology of the support function on the fan.
In the case for invertible sheaves on a complete toric variety, this is
a well-known fact.
In Section~$\ref{s3}$, we introduce the logarithmic differential forms
on a simplicial toric variety with an equivariant morphism to an affine
toric variety.
Under some assumption for the singularity of the simplicial toric
variety, that is caused by the characteristic of the base field, we
construct the Poincar\'{e} residue map for the sheaf of relative
logarithmic differential forms.
Using the Poincar\'{e} residue map, we prove the Bott vanishing theorem
and Euler exact sequence.
In Section~$\ref{s4}$, we consider hypersurfaces in a nonsingular toric
variety with a logarithmically smooth proper equivariant morphism to an
affine toric variety.
We define the Jacobian rings for hypersurfaces over the affine toric
variety, and prove the main result for describing the cohomology of the
sheaf of relative logarithmic differential forms by the Jacobian rings.
\section{Invertible sheaves on toric varieties}\label{s2}
First we introduce basic notation used in this paper, and then give some
properties of invertible sheaves on a toric variety with a proper
equivariant morphism to an affine toric variety.
We refer to \cite{f} and \cite{o} for terminology and basic facts in
toric geometry.\par
Let $N$ be a finitely generated free $\mathbf{Z}$-module.
We denote by $N_{\mathbf{R}}$ the $\mathbf{R}$-vector space
$\mathbf{R}\otimes_{\mathbf{Z}}N$, denote by $M$ the dual
$\mathbf{Z}$-module of $N$, and denote by
$\langle\ ,\ \rangle:M_{\mathbf{R}}\times N_{\mathbf{R}}
\rightarrow\mathbf{R}$
the canonical bilinear form.
Let $\sigma$ be a strongly convex rational polyhedral cone in
$N_{\mathbf{R}}$.
The dual cone of $\sigma$ is defined as
$$
\sigma^{\vee}=\{u\in M_{\mathbf{R}}\mid\langle u,v\rangle\geq0\
\text{for any}\ v\in\sigma\},
$$
and we denote by
$\mathbf{A}_{\sigma}=\Spec{k[M\cap\sigma^{\vee}]}$
the affine toric variety associated to $\sigma$ over an algebraically
closed field $k$.
Let $\Sigma$ be a finite fan of strongly convex rational polyhedral
cones in $N_{\mathbf{R}}$.
We denote by $|\Sigma|=\bigcup_{\sigma\in\Sigma}\sigma$ the support of
$\Sigma$, and denote by
$\mathbf{P}_{\Sigma}=\bigcup_{\sigma\in\Sigma}\mathbf{A}_{\sigma}$
the toric variety associated to $\Sigma$ over $k$.
Then the algebraic torus
$\mathbf{T}_{N}=\Spec{k[M]}$
is contained in $\mathbf{P}_{\Sigma}$ as an open subvariety, and
$\mathbf{T}_{N}$ acts on $\mathbf{P}_{\Sigma}$ as an extension of the
translation of $\mathbf{T}_{N}$.
For $u\in M$, the corresponding character
$$
\chi^{u}:\mathbf{T}_{N}\longrightarrow
\mathbf{G}_{\mathrm{m}}=\Spec{k[\mathbf{Z}]}
$$
is considered to be a rational function on $\mathbf{P}_{\Sigma}$.\par
We denote by $\Sigma(r)$ the set of all $r$-dimensional cones in
$\Sigma$.
For $\tau\in\Sigma(r)$, we define a free $\mathbf{Z}$-module by
$N_{\tau}=N/(N\cap\tau_{\mathbf{R}})$,
where $\tau_{\mathbf{R}}$ is the subspace of $N_{\mathbf{R}}$ generated
by $\tau$ over $\mathbf{R}$, and we define a set of subsets in
$N_{\tau,\mathbf{R}}$ by
$$
\Sigma_{\tau}=\{[\sigma]_{\tau}\mid\sigma\in\Sigma\
\text{and}\ \sigma\supset\tau\},
$$
where $[\sigma]_{\tau}$ is the image of $\sigma$ by the natural
homomorphism
$$
N_{\mathbf{R}}\longrightarrow N_{\tau,\mathbf{R}};\
v\longmapsto[v]_{\tau}.
$$
Then $\Sigma_{\tau}$ is a finite fan of strongly convex rational
polyhedral cones in $N_{\tau,\mathbf{R}}$, and the associated toric
variety $\mathbf{P}_{\Sigma_{\tau}}$ can be considered a
$\mathbf{T}_{N}$-invariant closed subvariety of codimension $r$ in
$\mathbf{P}_{\Sigma}$.
The closed immersion
$\iota_{\tau}:\mathbf{P}_{\Sigma_{\tau}}\rightarrow
\mathbf{P}_{\Sigma}$
is induced by
$$
\iota_{\tau}^{*}:k[M]\longrightarrow k[M_{\tau}];\
\chi^{u}\longmapsto
\begin{cases}
 \chi^{u}_{\tau}&(u\in M\cap\tau^{\perp}),\\
 0&(u\notin M\cap\tau^{\perp}),
\end{cases}
$$
where
$\chi^{u}_{\tau}:\mathbf{T}_{N_{\tau}}
\rightarrow\mathbf{G}_{\mathrm{m}}$
is the character corresponding to $u\in M_{\tau}=M\cap\tau^{\perp}$.\par
Let $B=\sum_{\rho\in\Sigma(1)}b_{\rho}\mathbf{P}_{\Sigma_{\rho}}$ be a
$\mathbf{T}_{N}$-invariant Weil divisor on $\mathbf{P}_{\Sigma}$.
We define a convex subset of $M_{\mathbf{R}}$ by
$$
\Delta_{B}=\{u\in M_{\mathbf{R}}\mid\langle u,v\rangle+b_{\rho}\geq0\
\text{for any}\ v\in|\Sigma|\}.
$$
Then there is a natural isomorphism
$$
H^{0}(\mathbf{P}_{\Sigma},\mathcal{O}_{\mathbf{P}_{\Sigma}}(B))
\simeq\bigoplus_{u\in M\cap\Delta_{B}}k\cdot\chi^{u}.
$$
Let $h:|\Sigma|\rightarrow\mathbf{R}$ be a $\Sigma$-linear support
function.
We define a $\mathbf{T}_{N}$-invariant Cartier divisor on
$\mathbf{P}_{\Sigma}$ by
$D_{h}=-\sum_{\rho\in\Sigma(1)}h(v_{\rho})\mathbf{P}_{\Sigma_{\rho}}$,
where $v_{\rho}$ is the generator of the monoid $\rho\cap N$.\par
Let $N'$ be another finitely generated free $\mathbf{Z}$-module, and let
$\pi_{*}:N\rightarrow N'$
be a homomorphism of $\mathbf{Z}$-modules.
We denote by
$\pi_{*\mathbf{R}}:N_{\mathbf{R}}\rightarrow
N'_{\mathbf{R}}$
the $\mathbf{R}$-homomorphism induced by $\pi_{*}$, and we denote by
$\pi^{*}:M'\rightarrow M$ the dual homomorphism of $\pi_{*}$.
Then $\pi^{*}$ induces a homomorphism of algebraic tori
$\pi_{0}:\mathbf{T}_{N}\rightarrow\mathbf{T}_{N'}$.
If a strongly convex rational polyhedral cone $\sigma'$ in
$N'_{\mathbf{R}}$ satisfies the condition
$|\Sigma|\subset\pi_{*\mathbf{R}}^{-1}(\sigma')$,
then $\pi^{*}$ induces an equivariant morphism of toric varieties
$\pi_{\Sigma,\sigma'}:\mathbf{P}_{\Sigma}
\rightarrow\mathbf{A}_{\sigma'}$,
which is an extension of $\pi_{0}$.
\begin{remark}
 The morphism $\pi_{\Sigma,\sigma'}$ is proper if and only if
 $|\Sigma|=\pi_{*\mathbf{R}}^{-1}(\sigma')$.
\end{remark}
\begin{remark}
 If $\pi_{\Sigma,\sigma'}$ is a proper morphism, then $|\Sigma|$ is a
 convex subset in $N_{\mathbf{R}}$.
 Conversely, if $|\Sigma|$ is a convex subset in $N_{\mathbf{R}}$, then
 we can find a free $\mathbf{Z}$-module $N'$, a surjective homomorphism
 $\pi_{*}:N\rightarrow N'$, and a strongly convex rational polyhedral
 cone $\sigma'$ in $N'_{\mathbf{R}}$ satisfying
 $|\Sigma|=\pi_{*\mathbf{R}}^{-1}(\sigma')$.
\end{remark}
\begin{theorem}\label{nef}
 Let $\mathbf{P}_{\Sigma}$ be a toric variety, let
 $\mathbf{A}_{\sigma'}$ be an affine toric variety, and let
 $\pi:\mathbf{P}_{\Sigma}\rightarrow\mathbf{A}_{\sigma'}$ be a
 proper equivariant morphism.
 For a $\Sigma$-linear support function $h$, the following conditions
 are equivalent:
 \renewcommand{\labelenumi}{$(\arabic{enumi})$\hspace*{5pt}}
 \begin{enumerate}
  \item $\mathcal{O}_{\mathbf{P}_{\Sigma}}(D_{h})$ is generated by
	global sections;
  \item $\mathcal{O}_{\mathbf{P}_{\Sigma}}(D_{h})$ is $\pi$-nef;
  \item $h$ is upper convex, i.e. for any $v_{1},v_{2}\in|\Sigma|$,
	$$
	h(v_{1}+v_{2})\geq h(v_{1})+h(v_{2}).
	$$
 \end{enumerate}
\end{theorem}
\begin{proof}
 $(1)\Rightarrow(2)$.\hspace*{5pt}
 Let $C$ be a complete integral curve in a fiber of $\pi$.
 We denote by $\widetilde{C}$ the normalization of $C$, and denote by
 $\iota$ the morphism from $\widetilde{C}$ to $\mathbf{P}_{\Sigma}$.
 Since $\mathcal{O}_{\mathbf{P}_{\Sigma}}(D_{h})$ is generated by global
 sections, $\iota^{*}\mathcal{O}_{\mathbf{P}_{\Sigma}}(D_{h})$ has a
 non-zero global section, so the intersection number is
 $$
 (D_{h}.\ C)
 =\deg{\iota^{*}\mathcal{O}_{\mathbf{P}_{\Sigma}}(D_{h})}
 \geq0.
 $$
 $(2)\Rightarrow(3)$.\hspace*{5pt}
 The morphism $\pi$ is defined as $\pi=\pi_{\Sigma,\sigma'}$ by a
 homomorphism
 $\pi_{*}:N\rightarrow N'$.
 We denote by $N''$ the image of $\pi_{*}$, and define a strongly convex
 rational polyhedral cone in $N''_{\mathbf{R}}$ by
 $\sigma''=\sigma'\cap N''_{\mathbf{R}}$.
 Then $\pi$ is factored to a proper surjective morphism
 $\pi_{\Sigma,\sigma''}:\mathbf{P}_{\Sigma}
 \rightarrow\mathbf{A}_{\sigma''}$,
 and a finite morphism
 $\mathbf{A}_{\sigma''}\rightarrow\mathbf{A}_{\sigma'}$.
 We denote by $s$ the dimension of the convex subset
 $$
 |\Sigma|=\pi_{*\mathbf{R}}^{-1}(\sigma')
 =\pi_{*\mathbf{R}}^{-1}(\sigma'')\subset N_{\mathbf{R}},
 $$
 and define a subset of $\Sigma(s-1)$ by
 $$
 \Sigma(\sigma'',s-1)=\{\tau\in\Sigma(s-1)\mid
 \pi_{*\mathbf{R}}(\tau)\cap\mathrm{Int}(\sigma'')\neq\emptyset\},
 $$
 where $\mathrm{Int}(\sigma'')$ denotes the relative interior of
 $\sigma''$.
 For $\tau\in\Sigma(\sigma'',s-1)$, there exist exactly two cones
 $\sigma_{+},\sigma_{-}\in\Sigma(s)$ containing $\tau$.
 Then
 $\pi_{\Sigma,\sigma''}\circ\iota_{\tau}:
 \mathbf{P}_{\Sigma_{\tau}}\rightarrow\mathbf{A}_{\sigma''}$
 gives a $\mathbf{P}^{1}$-bundle over the closed
 $\mathbf{T}_{N''}$-orbit
 $\mathbf{T}_{N''_{\sigma''}}=\Spec{k[M''\cap\sigma''^{\perp}]}$
 in $\mathbf{A}_{\sigma''}$, and
 $\mathbf{P}_{\Sigma_{\sigma_{+}}}$ and
 $\mathbf{P}_{\Sigma_{\sigma_{-}}}$ become the $0$-section and the
 $\infty$-section of the $\mathbf{P}^{1}$-bundle.
 Let $0\in\mathbf{T}_{N''_{\sigma''}}(k)$ be a $k$-rational point.
 We denote by $\mathbf{P}_{\Sigma_{\tau},0}$ the fiber of
 $\pi_{\Sigma,\sigma''}\circ\iota_{\tau}$ at $0$.
 Then $\mathbf{P}_{\Sigma_{\tau},0}$ is a
 nonsingular rational curve in a fiber of $\pi$.
 Since $\rho_{-}=[\sigma_{-}]_{\tau}$ is a $1$-dimensional cone in
 $N_{\tau,\mathbf{R}}$, for $v\in\sigma_{-}\smallsetminus\tau$, a
 positive real number $a_{v}$ is determined by
 $$
 [v]_{\tau}=a_{v}v_{\rho_{-}}
 \in\rho_{-}\subset N_{\tau,\mathbf{R}}
 =N_{\mathbf{R}}/\tau_{\mathbf{R}}.
 $$
 For $\sigma\in\Sigma(s)$, we denote by $u_{\sigma}\in M$ a linear
 function which coincides with $h$ on $\sigma$.
 Then we have an equation
 $$
 (D_{h}.\ \mathbf{P}_{\Sigma_{\tau},0})
 =\frac{1}{a_{v}}(\langle u_{\sigma_{+}},v\rangle-h(v))
 $$
 for any $v\in\sigma_{-}\smallsetminus\tau$.
 Since $\mathcal{O}_{\mathbf{P}_{\Sigma}}(D_{h})$ is $\pi$-nef, we have
 \begin{equation}\label{pnef}
  \langle u_{\sigma_{+}},v\rangle\geq h(v)
 \end{equation}
 for any $v\in\sigma_{-}$.\par
 Next we prove $(\ref{pnef})$ without restriction about $\sigma_{+}$ and
 $v$.
 For $\sigma\in\Sigma(s)$ and $v\in|\Sigma|$, there exists a vector
 $w\in\mathrm{Int}(\sigma)$ such that
 $$
 L(w,v)=\{(1-t)w+tv\in N_{\mathbf{R}}\mid0<t<1\}
 $$
 does not intersect with $\bigcup_{\tau\in\Sigma(s-2)}\tau$.
 Then $\bigcup_{\tau\in\Sigma(\sigma'',s-1)}\tau$ divides
 $L(w,v)$ into finite pieces, and we define vectors
 $w_{0},\ldots,w_{l}\in N_{\mathbf{R}}$ by
 $$
 L(w,v)\cap\bigcup_{\tau\in\Sigma(\sigma'',s-1)}\tau
 =\{w_{1},\ldots,w_{l-1}\}
 $$
 and
 $$
 w_{i}=(1-t_{i})w+t_{i}v
 \hspace*{20pt}
 (0=t_{0}<t_{1}<\cdots<t_{l-1}<t_{l}=1).
 $$
 For $1\leq i\leq l$, there exists a unique cone
 $\sigma_{i}\in\Sigma(s)$ such that $w_{i-1},w_{i}\in\sigma_{i}$.
 We note that $w=w_{0}$, $v=w_{l}$ and $\sigma=\sigma_{1}$.
 Since
 $$
 h(w_{i})=\langle u_{\sigma_{i}},w_{i}\rangle
 = \langle u_{\sigma_{i+1}},w_{i}\rangle
 \hspace*{20pt}
 (1\leq i\leq l-1),
 $$
 we have
 \begin{align*}
  \langle u_{\sigma},v\rangle
  &=\langle u_{\sigma_{l}},v\rangle+\sum_{i=1}^{l-1}
  \langle u_{\sigma_{i}}-u_{\sigma_{i+1}},v\rangle\\
  &=h(v)+\sum_{i=1}^{l-1}
  \langle u_{\sigma_{i}}-u_{\sigma_{i+1}},
  \frac{1-t_{i}}{t_{i+1}-t_{i}}w_{i+1}
  +\frac{t_{i+1}-1}{t_{i+1}-t_{i}},w_{i}\rangle\\
  &=h(v)+\sum_{i=1}^{l-1}
  \frac{1-t_{i}}{t_{i+1}-t_{i}}
  (\langle u_{\sigma_{i}},w_{i+1}\rangle-h(w_{i+1})).
 \end{align*}
 By using $(\ref{pnef})$ for $\sigma_{+}=\sigma_{i}$ and
 $\sigma_{-}=\sigma_{i+1}$, we have
 $$
 \langle u_{\sigma_{i}},w_{i+1}\rangle\geq h(w_{i+1})
 \hspace*{20pt}
 (1\leq i\leq l-1),
 $$
 hence
 \begin{equation}\label{uconv}
  \langle u_{\sigma},v\rangle
   \geq h(v)
 \end{equation}
 for any $\sigma\in\Sigma(s)$ and $v\in|\Sigma|$.\par
 For $v_{1},v_{2}\in|\Sigma|$, there is a cone $\sigma\in\Sigma(s)$ such
 that $v_{1}+v_{2}\in\sigma$.
 By $(\ref{uconv})$,
 $$
 h(v_{1}+v_{2})=\langle u_{\sigma},v_{1}+v_{2}\rangle
 =\langle u_{\sigma},v_{1}\rangle+
 \langle u_{\sigma},v_{2}\rangle
 \geq h(v_{1})+h(v_{2}).
 $$
 $(3)\Rightarrow(1)$.\hspace*{5pt}
 For $\sigma\in\Sigma(s)$ and $v\in|\Sigma|$, there exists $w\in\sigma$
 such that $v+w$ is contained in $\sigma$.
 Because $h$ is upper convex, we have
 $$
 \langle u_{\sigma},v\rangle
 =\langle u_{\sigma},v+w\rangle-\langle u_{\sigma},w\rangle
 =h(v+w)-h(w)\geq h(v),
 $$
 so $\chi^{u_{\sigma}}$ gives a global section of
 $\mathcal{O}_{\mathbf{P}_{\Sigma}}(D_{h})$.
 Since
 $\mathbf{P}_{\Sigma}=\bigcup_{\sigma\in\Sigma(s)}\mathbf{A}_{\sigma}$,
 and $\chi^{u_{\sigma}}$ generates
 $\varGamma(\mathbf{A}_{\sigma},
 \mathcal{O}_{\mathbf{P}_{\Sigma}}(D_{h}))$
 over $k[M\cap\sigma^{\vee}]$, the invertible sheaf
 $\mathcal{O}_{\mathbf{P}_{\Sigma}}(D_{h})$ is generated by global
 sections.
\end{proof}
\begin{theorem}\label{ample}
 Let $\mathbf{P}_{\Sigma}$ be a toric variety, let
 $\mathbf{A}_{\sigma'}$ be an affine toric variety, and let
 $\pi:\mathbf{P}_{\Sigma}\rightarrow\mathbf{A}_{\sigma'}$ be a
 proper equivariant morphism.
 For a $\Sigma$-linear support function $h$, the following conditions
 are equivalent:
 \renewcommand{\labelenumi}{$(\arabic{enumi})$\hspace*{5pt}}
 \begin{enumerate}
  \item $\mathcal{O}_{\mathbf{P}_{\Sigma}}(D_{h})$ is ample;
  \item $\mathcal{O}_{\mathbf{P}_{\Sigma}}(D_{h})$ is $\pi$-ample;
  \item $h$ is strictly upper convex, i.e. for any
	$v_{1},v_{2}\in|\Sigma|$,
	$$
	h(v_{1}+v_{2})\geq h(v_{1})+h(v_{2}),
	$$
	and equality holds if and only if there exists a cone
	$\sigma\in\Sigma$ such that $v_{1},v_{2}\in\sigma$.
 \end{enumerate}
\end{theorem}
\begin{proof}
 The equivalence of $(1)$ and $(2)$ is well-known for any proper
 morphism to an affine scheme.\\
 $(2)\Rightarrow(3)$.\hspace*{5pt}
 In the proof of Theorem~$\ref{nef}$ $(2)\Rightarrow(3)$, if
 $\mathcal{O}_{\mathbf{P}_{\Sigma}}(D_{h})$ is $\pi$-ample, then
 $(D_{h}.\ \mathbf{P}_{\Sigma_{\tau},0})>0$,
 hence
 $\langle u_{\sigma_{+}},v\rangle>h(v)$
 for any $v\in\sigma_{-}\smallsetminus\tau$.
 This implies that
 $\langle u_{\sigma},v\rangle\geq h(v)$
 for any $\sigma\in\Sigma(s)$ and $v\in|\Sigma|$,
 and equality holds if and only if $v\in\sigma$.\par
 For $v_{1},v_{2}\in|\Sigma|$, there is a cone $\sigma\in\Sigma(s)$ such
 that $v_{1}+v_{2}\in\sigma$.
 If
 $$
 h(v_{1}+v_{2})=h(v_{1})+h(v_{2}),
 $$
 then
 $$
 (\langle u_{\sigma},v_{1}\rangle-h(v_{1}))+
 (\langle u_{\sigma},v_{2}\rangle-h(v_{2}))=0.
 $$
 So we have
 $\langle u_{\sigma},v_{1}\rangle=h(v_{1})$
 and
 $\langle u_{\sigma},v_{2}\rangle=h(v_{2})$,
 which mean $v_{1},v_{2}\in\sigma$.\\
 $(3)\Rightarrow(1)$.\hspace*{5pt}
 Let $\mathcal{E}$ be a coherent
 $\mathcal{O}_{\mathbf{P}_{\Sigma}}$-module.
 We show that there exists a positive integer $m_{0}$ such that
 $\mathcal{E}\otimes_{\mathcal{O}_{\mathbf{P}_{\Sigma}}}
 \mathcal{O}_{\mathbf{P}_{\Sigma}}(mD_{h})$ is
 generated by global sections for any integer $m\geq m_{0}$.
 We may assume that
 $\mathcal{E}=\mathcal{O}_{\mathbf{P}_{\Sigma}}(B)$
 for a $\mathbf{T}_{N}$-invariant Weil divisor
 $B=\sum_{\rho\in\Sigma(1)}b_{\rho}\mathbf{P}_{\Sigma_{\rho}}$,
 because by \cite{m},
 there exists a surjective homomorphism
 $\bigoplus_{i=1}^{r}\mathcal{O}_{\mathbf{P}_{\Sigma}}(B_{i})
 \rightarrow\mathcal{E}$
 for some $\mathbf{T}_{N}$-invariant Weil divisors
 $B_{1},\ldots,B_{r}$.\par
 For $\sigma\in\Sigma(s)$, we fix vectors
 $u_{\sigma,1},\ldots,u_{\sigma,c_{\sigma}}\in M$
 such that $\chi^{u_{\sigma,1}},\ldots,\chi^{u_{\sigma,c_{\sigma}}}$
 generate
 $\varGamma(\mathbf{A}_{\sigma},
 \mathcal{O}_{\mathbf{P}_{\Sigma}}(B))$
 over $k[M\cap\sigma^{\vee}]$.
 Since $h$ is strictly upper convex, for $\sigma\in\Sigma(s)$ and
 $\rho\in\Sigma(1)$, we have
 $\langle u_{\sigma},v_{\rho}\rangle\geq h(v_{\rho})$,
 and if $v_{\rho}\notin\sigma$, then
 $\langle u_{\sigma},v_{\rho}\rangle>h(v_{\rho})$.
 And if $v_{\rho}\in\sigma$, then
 $\langle u_{\sigma,i},v_{\rho}\rangle+b_{\rho}\geq0$.
 Hence there exists a positive integer $m_{0}$ such that
 for any $\sigma\in\Sigma(s)$, for any $1\leq i\leq c_{\sigma}$ and for
 any $\rho\in\Sigma(1)$,
 $$
 m_{0}(\langle u_{\sigma},v_{\rho}\rangle-h(v_{\rho}))\geq
 -b_{\rho}-\langle u_{\sigma,i},v_{\rho}\rangle.
 $$
 Then we have
 $$
 \langle u_{\sigma,i}+mu_{\sigma},v_{\rho}\rangle
 \geq -b_{\rho}+mh(v_{\rho})
 $$
 for any $m\geq m_{0}$, and this means that
 $\chi^{u_{\sigma,i}+mu_{\sigma}}$ is a global section of the coherent
 sheaf $\mathcal{O}_{\mathbf{P}_{\Sigma}}(B+mD_{h})$.
 Since
 $\chi^{u_{\sigma,1}+mu_{\sigma}},\ldots,
 \chi^{u_{\sigma,c_{\sigma}}+mu_{\sigma}}$
 generate
 $\varGamma(\mathbf{A}_{\sigma},
 \mathcal{O}_{\mathbf{P}_{\Sigma}}(B+mD_{h}))$
 over $k[M\cap\sigma^{\vee}]$, the coherent sheaf
 $\mathcal{O}_{\mathbf{P}_{\Sigma}}(B+mD_{h})$ is
 generated by global sections.
\end{proof}
\section{Log differential forms on toric varieties}\label{s3}
We introduce the sheaf of relative logarithmic
differential forms on a toric variety with an equivariant morphism to an
affine toric variety.\par
Let $\mathbf{P}=\mathbf{P}_{\Sigma}$ be a toric variety, let
$\mathbf{A}=\mathbf{A}_{\sigma'}$ be an affine toric variety, and let
$\pi=\pi_{\Sigma,\sigma'}:\mathbf{P}\rightarrow\mathbf{A}$ be an
equivariant morphism, which is given by a homomorphism
$\pi_{*}:N\rightarrow N'$ with
$|\Sigma|\subset\pi_{*\mathbf{R}}^{-1}(\sigma')$.
For $\tau\in\Sigma(r)$, we denote by $\mathbf{P}_{\tau}$ the
corresponding $\mathbf{T}_{N}$-invariant subvariety
$\mathbf{P}_{\Sigma_{\tau}}$.
We define a subset of $\Sigma(r)$ by
$$
\Sigma_{\pi}(r)=\{\tau\in\Sigma(r)
\mid\tau\subset\pi_{*\mathbf{R}}^{-1}(0)\}.
$$
If $\tau\in\Sigma_{\pi}(r)$, then
$\pi_{\tau}=\pi\circ\iota_{\tau}:
\mathbf{P}_{\tau}\rightarrow\mathbf{A}$
is an equivariant morphism induced by
$$
N_{\tau}\longrightarrow N';\
[v]_{\tau}\longmapsto\pi_{*}(v).
$$
We denote by $\Sigma^{\mathrm{reg}}$ the set of all nonsingular cones in
$\Sigma$, and denote by
$j:\mathbf{P}^{\mathrm{reg}}=\mathbf{P}_{\Sigma^{\mathrm{reg}}}
\rightarrow\mathbf{P}$
the natural open immersion.
We define $\mathbf{T}_{N}$-invariant divisors on $\mathbf{P}$ by
\begin{equation}\label{div}
 \begin{cases}
  D=
  \begin{displaystyle}
   \sum_{\rho\in\Sigma_{\pi}(1)}\mathbf{P}_{\rho},
  \end{displaystyle}\\
  E=
  \begin{displaystyle}
   \sum_{\rho\in\Sigma(1)\smallsetminus\Sigma_{\pi}(1)}
   \mathbf{P}_{\rho},
  \end{displaystyle}
 \end{cases}
\end{equation}
which are divisors with normal crossing on
$\mathbf{P}^{\mathrm{reg}}$.\par
Let $\sigma$ be a cone in $\Sigma^{\mathrm{reg}}(c)$, and let
$(v_{1},\ldots,v_{d})$ be a $\mathbf{Z}$-basis of $N$ with
\begin{equation}\label{local}
 \begin{cases}
  \sigma=\mathbf{R}_{\geq0}v_{1}+\cdots+\mathbf{R}_{\geq0}v_{c},\\
  v_{1},\ldots,v_{l}\in\pi_{*}^{-1}(0),\
  v_{l+1},\ldots,v_{c}\notin\pi_{*}^{-1}(0).
 \end{cases}
\end{equation}
Then we have
\begin{align*}
 &\mathbf{A}_{\sigma}=\Spec{k[M\cap\sigma^{\vee}]}
 =\Spec{k[x_{1},\ldots,x_{c},
 x_{c+1}^{\pm1},\ldots,x_{d}^{\pm1}]},\\
 &\mathbf{A}_{\sigma}\cap D=\Spec{k[x_{1},\ldots,x_{c},
 x_{c+1}^{\pm1},\ldots,x_{d}^{\pm1}]/(x_{1}\cdots x_{l})},\\
 &\mathbf{A}_{\sigma}\cap E=\Spec{k[x_{1},\ldots,x_{c},
 x_{c+1}^{\pm1},\ldots,x_{d}^{\pm1}]/(x_{l+1}\cdots x_{c})},
\end{align*}
where $x_{i}=\chi^{u_{i}}$ for the dual basis $(u_{1},\ldots,u_{d})$ of
$(v_{1},\ldots,v_{d})$.
We define a free $k[M\cap\sigma^{\vee}]$-module by
$$
\omega_{\mathbf{A}_{\sigma}}^{1}
=\varOmega_{\mathbf{A}_{\sigma}}^{1}
(\log{E})
=\bigoplus_{j=1}^{l}k[M\cap\sigma^{\vee}]dx_{j}
\oplus\bigoplus_{j=l+1}^{d}k[M\cap\sigma^{\vee}]
\frac{dx_{j}}{x_{j}},
$$
which is naturally contained in the free $k[M\cap\sigma^{\vee}]$-module
$$
\omega_{\mathbf{A}_{\sigma}}^{1}(\log{D})
=\varOmega_{\mathbf{A}_{\sigma}}^{1}
(\log{D\cup E})
=\bigoplus_{j=1}^{d}k[M\cap\sigma^{\vee}]
\frac{dx_{j}}{x_{j}}.
$$
We denote by
$\omega_{\mathbf{P}^{\mathrm{reg}}}^{1}\subset
\omega_{\mathbf{P}^{\mathrm{reg}}}^{1}(\log{D})$
the sheaves of $\mathcal{O}_{\mathbf{P}^{\mathrm{reg}}}$-modules defined
by $k[M\cap\sigma^{\vee}]$-modules
$\omega_{\mathbf{A}_{\sigma}}^{1}\subset
\omega_{\mathbf{A}_{\sigma}}^{1}(\log{D})$
for $\sigma\in\Sigma^{\mathrm{reg}}$.
Then we have an isomorphism
$\mathcal{O}_{\mathbf{P}^{\mathrm{reg}}}\otimes_{\mathbf{Z}}M
\simeq
\omega_{\mathbf{P}^{\mathrm{reg}}}^{1}(\log{D})$
by
$$
k[M\cap\sigma^{\vee}]\otimes_{\mathbf{Z}}M
\simeq
\omega_{\mathbf{A}_{\sigma}}^{1}(\log{D});\
1\otimes u_{j}\longleftrightarrow\frac{dx_{j}}{x_{j}}.
$$
We denote by
$\omega_{\mathbf{P}^{\mathrm{reg}}/\mathbf{A}}^{1}(\log{D})$
the cokernel of the homomorphism
$$
\pi^{*}_{\mathcal{O}_{\mathbf{P}^{\mathrm{reg}}}}:
\mathcal{O}_{\mathbf{P}^{\mathrm{reg}}}\otimes_{\mathbf{Z}}M'
\longrightarrow
\mathcal{O}_{\mathbf{P}^{\mathrm{reg}}}\otimes_{\mathbf{Z}}M
\simeq\omega_{\mathbf{P}^{\mathrm{reg}}}^{1}(\log{D}).
$$
Since the image of $\pi^{*}_{\mathcal{O}_{\mathbf{P}^{\mathrm{reg}}}}$
is contained in $\omega_{\mathbf{P}^{\mathrm{reg}}}^{1}$, we denote by
$\omega_{\mathbf{P}^{\mathrm{reg}}/\mathbf{A}}^{1}$ the cokernel of
$\pi^{*}_{\mathcal{O}_{\mathbf{P}^{\mathrm{reg}}}}:
\mathcal{O}_{\mathbf{P}^{\mathrm{reg}}}\otimes_{\mathbf{Z}}M'
\rightarrow\omega_{\mathbf{P}^{\mathrm{reg}}}^{1}$.
We define a coherent sheaf on $\mathbf{P}$ by
$\widetilde{\omega}_{\mathbf{P}/\mathbf{A}}^{p}
=j_{*}(\bigwedge^{p}
\omega_{\mathbf{P}^{\mathrm{reg}}/\mathbf{A}}^{1})$,
which is a submodule of the free $\mathcal{O}_{\mathbf{P}}$-module
$$
\omega_{\mathbf{P}/\mathbf{A}}^{p}(\log{D})=
j_{*}\omega_{\mathbf{P}^{\mathrm{reg}}/\mathbf{A}}^{p}(\log{D})
\simeq\mathcal{O}_{\mathbf{P}}\otimes_{k}\bigwedge^{p}M_{\pi,k},
$$
where $M_{\pi,k}$ denotes the cokernel of
$\pi^{*}_{k}:k\otimes_{\mathbf{Z}}M'\rightarrow
k\otimes_{\mathbf{Z}}M$.
The sheaf $\widetilde{\omega}_{\mathbf{P}/\mathbf{A}}^{p}$ is the sheaf
of relative logarithmic differential $p$-forms of Zariski.
In the paper \cite{bc}, it is simply denoted by
$\varOmega_{\mathbf{P}}^{p}$ for the case $\mathbf{A}=\Spec{k}$.
\begin{remark}
 The shaves $\omega_{\mathbf{P}^{\mathrm{reg}}/\mathbf{A}}^{1}$ and
 $\omega_{\mathbf{P}/\mathbf{A}}^{1}(\log{D})$ are interpreted as
 sheaves of relative logarithmic differential forms in the sense of
 logarithmic geometry \cite{k}.
 If we consider $\mathbf{P}^{\mathrm{reg}}$ with the logarithmic
 structure $\mathcal{N}_{\mathbf{P}^{\mathrm{reg}}}^{E}$ defined by the
 divisor $E$ with normal crossing, and consider $\mathbf{A}$ with the
 canonical logarithmic structure
 $\mathcal{N}_{\mathbf{A}}^{\mathrm{can}}$ as a toric variety, then the
 sheaf $\omega_{\mathbf{P}^{\mathrm{reg}}/\mathbf{A}}^{1}$ is the sheaf
 of differential forms on
 $(\mathbf{P}^{\mathrm{reg}},
 \mathcal{N}_{\mathbf{P}^{\mathrm{reg}}}^{E})$
 over $(\mathbf{A},\mathcal{N}_{\mathbf{A}}^{\mathrm{can}})$.
 If we consider $\mathbf{P}$ with the canonical logarithmic structure
 $\mathcal{N}_{\mathbf{P}}^{\mathrm{can}}$ as a toric variety, then
 the sheaf $\omega_{\mathbf{P}/\mathbf{A}}^{1}(\log{D})$ is the
 sheaf of differential forms on
 $(\mathbf{P},\mathcal{N}_{\mathbf{P}}^{\mathrm{can}})$ over
 $(\mathbf{A},\mathcal{N}_{\mathbf{A}}^{\mathrm{can}})$.
\end{remark}
We describe sections of the sheaf
$\widetilde{\omega}_{\mathbf{P}/\mathbf{A}}^{p}$
by the following explicit way, that is the idea in \cite{da}.
For $u\in M$ and $\sigma\in\Sigma$, we define a subset of
$\Sigma_{\pi}(1)$ by
$$
\sigma_{\pi,u=0}(1)=\{\rho\in\Sigma_{\pi}(1)\mid
\langle u,v_{\rho}\rangle=0\ \text{and}\ \rho\subset\sigma\},
$$
and define a $k$-subspace of $M_{\pi,k}$ by
$$
H_{\sigma,u}=\{w\in M_{\pi,k}\mid
\langle w,1\otimes v_{\rho}\rangle=0\
\text{for any}\ \rho\in\sigma_{\pi,u=0}(1)\}.
$$
\begin{lemma}\label{li}
 For $\sigma\in\Sigma$, there is a natural isomorphism
 $$
 \varGamma(\mathbf{A}_{\sigma},
 \widetilde{\omega}_{\mathbf{P}/\mathbf{A}}^{p})
 \simeq\bigoplus_{u\in M\cap\sigma^{\vee}}\bigwedge^{p}
 H_{\sigma,u}\cdot\chi^{u}.
 $$
\end{lemma}
\begin{proof}
 The isomorphism is induced from the natural isomorphism
 $$
 \varGamma(\mathbf{A}_{\sigma},
 \omega_{\mathbf{P}/\mathbf{A}}^{p}(\log{D}))
 \simeq\varGamma(\mathbf{A}_{\sigma},
 \mathcal{O}_{\mathbf{P}})\otimes_{k}\bigwedge^{p}M_{\pi,k}
 \simeq\bigoplus_{u\in M\cap\sigma^{\vee}}
 \bigwedge^{p}M_{\pi,k}\cdot\chi^{u}.
 $$
 First we assume that $\sigma\in\Sigma^{\mathrm{reg}}$.
 We use a $\mathbf{Z}$-basis $(v_{1},\ldots,v_{d})$ of $N$ satisfying
 $(\ref{local})$ and the dual basis $(u_{1},\ldots,u_{d})$ of $M$.
 There exists a $k$-basis $(w_{1},\ldots,w_{n})$ of $M_{\pi,k}$ such
 that for $1\leq i\leq l$, the vector $w_{i}$ is the image of
 $1\otimes u_{i}$ by the natural homomorphism
 $k\otimes_{\mathbf{Z}}M\rightarrow M_{\pi,k}$, and
 $\varGamma(\mathbf{A}_{\sigma},\omega_{\mathbf{P}^{\mathrm{reg}}
 /\mathbf{A}}^{1})$
 corresponds to the $k[M\cap\sigma^{\vee}]$-submodule
 $H_{\sigma}\subset
 \bigoplus_{u\in M\cap\sigma^{\vee}}M_{\pi,k}\cdot\chi^{u}$
 generated by
 $$
 w_{1}\chi^{u_{1}},\ldots,w_{l}\chi^{u_{l}},
 w_{l+1}\chi^{0},\ldots,w_{n}\chi^{0}.
 $$
 For a subset $I=\{i_{1},\ldots,i_{p}\}\subset\{1,\ldots,n\}$
 with $i_{1}<\cdots<i_{p}$, we define an element in $M$ by
 $u_{I}=\sum_{i\in I\cap\{1,\ldots,l\}}u_{i}$,
 and define a vector in $\bigwedge^{p}M_{\pi,k}$ by
 \begin{equation}\label{mi}
  w_{I}=w_{i_{1}}\wedge\cdots\wedge w_{i_{p}}.
 \end{equation}
 Then $(w_{I};\ |I|=p)$ is a $k$-basis of $\bigwedge^{p}M_{\pi,k}$, and
 $(w_{I}\chi^{u_{I}};\ |I|=p)$ is a $k[M\cap\sigma^{\vee}]$-basis of
 $\bigwedge^{p}H_{\sigma}$.
 For $u\in M\cap\sigma^{\vee}$, we define a subset of $\{1,\ldots,l\}$
 by
 $$
 J_{u}=\{j\in\{1,\ldots,l\}\mid\langle u,v_{j}\rangle=0\}.
 $$
 Then $(w_{I};\ |I|=p,\ I\cap J_{u}=\emptyset)$ is a $k$-basis of
 $\bigwedge^{p}H_{\sigma,u}$.
 If
 $$
 \omega=
 \sum_{u\in M\cap\sigma^{\vee}}\sum_{|I|=p}a_{u,I}w_{I}\chi^{u}\in
 \bigoplus_{u\in M\cap\sigma^{\vee}}\bigwedge^{p}M_{\pi,k}\cdot\chi^{u}
 \hspace*{20pt}(a_{u,I}\in k)
 $$
 is contained in $\bigwedge^{p}H_{\sigma}$,
 then
 $\omega=
 \sum_{I=|p|}(\sum_{u\in M\cap\sigma^{\vee}}b_{I,u}\chi^{u})
 w_{I}\chi^{u_{I}}$
 for some $b_{I,u}\in k$, so we have
 $$
 a_{u,I}=
 \begin{cases}
  b_{I,u-u_{I}}&(u-u_{I}\in\sigma^{\vee})\\
  0&(u-u_{I}\notin\sigma^{\vee}).
 \end{cases}
 $$
 If $u-u_{I}\in\sigma^{\vee}$, then $I\cap J_{u}=\emptyset$, hence
 $\sum_{|I|=p}a_{u,I}w_{I}\in\bigwedge^{p}H_{\sigma,u}$.
 Conversely, if $I\cap J_{u}=\emptyset$, then we have
 $u-u_{I}\in\sigma^{\vee}$, hence
 $$
 w_{I}\chi^{u}=\chi^{u-u_{I}}(w_{I}\chi^{u_{I}})
 \in\bigwedge^{p}H_{\sigma}.
 $$
 So $\bigwedge^{p}H_{\sigma}$ coincides with
 $\bigoplus_{u\in M\cap\sigma^{\vee}}\bigwedge^{p}
 H_{\sigma,u}\cdot\chi^{u}$
 for $\sigma\in\Sigma^{\mathrm{reg}}$.\par
 If $\sigma\in\Sigma$ is not a nonsingular cone, then we have
 $$
 \varGamma(\mathbf{A}_{\sigma},
 \widetilde{\omega}_{\mathbf{P}/\mathbf{A}}^{p})
 =\varGamma(\mathbf{A}_{\sigma}\cap\mathbf{P}^{\mathrm{reg}},
 \omega_{\mathbf{P}^{\mathrm{reg}}/\mathbf{A}}^{p})
 =\bigcap_{\tau\in\sigma^{\mathrm{reg}}}
 \varGamma(\mathbf{A}_{\tau},\omega_{\mathbf{P}^{\mathrm{reg}}
 /\mathbf{A}}^{p}),
 $$
 where $\sigma^{\mathrm{reg}}$ denotes the set of all nonsingular faces
 of $\sigma$.
 Since
 $\sigma^{\vee}=\bigcap_{\tau\in\sigma^{\mathrm{reg}}}\tau^{\vee}$,
 and
 $\bigwedge^{p}H_{\sigma,u}
 =\bigcap_{\tau\in\sigma^{\mathrm{reg}}}
 \bigwedge^{p}H_{\tau,u}$
 for $u\in\sigma^{\vee}$, the isomorphism in Lemma~$\ref{li}$ is proved
 for any $\sigma\in\Sigma$.
\end{proof}
Next we define the Poincar\'{e} residue map for the sheaf
$\widetilde{\omega}_{\mathbf{P}/\mathbf{A}}^{p}$.
For this purpose, we need some assumptions for the fan $\Sigma$.
For $\sigma\in\Sigma$, we define a positive integer by
$$
m(\sigma)=
\ord{(N\cap\sigma_{\mathbf{R}}/
\sum_{\rho\in\sigma(1)}\mathbf{Z}v_{\rho})},
$$
and define a positive integer by
$m(\Sigma)=\prod_{\sigma\in\Sigma}m(\sigma)$.
We assume that the fan $\Sigma$ is simplicial, and assume that
$m(\Sigma)$ is prime to the characteristic of $k$.
For $\tau\in\Sigma_{\pi}(r)$, we define a $k$-homomorphism by
$$
\phi_{\tau}:\bigwedge^{r}M_{\pi,k}\longrightarrow k;\
\eta_{1}\wedge\cdots\wedge\eta_{r}\longmapsto
\det{(\langle\eta_{i},1\otimes v_{j}\rangle)
_{1\leq i\leq r,1\leq j\leq r}},
$$
where we fix an ordering
$$
\{v_{1},\ldots,v_{r}\}=\{v_{\rho}\in N\mid\rho\in\Sigma_{\pi}(1)\
\text{and}\ \rho\subset\tau\}.
$$
The homomorphism $\phi_{\tau}$ is depend on the ordering, but it is
determined modulo $\pm1$.
For $\sigma\in\Sigma$ and $\tau\in\Sigma_{\pi}(r)$ with
$\tau\subset\sigma$, we define a homomorphism by
\begin{align*}
 \Phi_{\sigma,\tau}:
 \bigwedge^{r}M_{\pi,k}&\otimes_{k}
 \bigoplus_{u\in M\cap\sigma^{\vee}}\bigwedge^{p-r}H_{\sigma,u}
 \cdot\chi^{u}
 &\longrightarrow&
 \bigoplus_{u\in M\cap\sigma^{\vee}\cap\tau^{\perp}}
 \bigwedge^{p-r}H_{\sigma,u}\cdot\chi^{u}_{\tau};\\
 &\eta\otimes w\chi^{u}
 &\longmapsto&
 \begin{cases}
  \phi_{\tau}(\eta)w\chi^{u}_{\tau}&(u\in\tau^{\perp}),\\
  0&(u\notin\tau^{\perp}).
 \end{cases}
\end{align*}
Since $m(\Sigma)$ is prime to the characteristic of $k$, for
$u\in M\cap\sigma^{\vee}\cap\tau^{\perp}
\simeq M_{\tau}\cap[\sigma]_{\tau}^{\vee}$,
there is a natural identification
$$
H_{\sigma,u}\simeq
\{w\in(M_{\tau})_{\pi_{\tau},k}\mid\langle w,1\otimes v_{\rho}\rangle=0\
\text{for any}\ \rho\in([\sigma]_{\tau})_{\pi_{\tau},u=0}(1)\},
$$
so $\Phi_{\sigma,\tau}$ gives a homomorphism
$$
\varGamma(\mathbf{A}_{\sigma},
\omega_{\mathbf{P}/\mathbf{A}}^{r}(\log{D})
\otimes_{\mathcal{O}_{\mathbf{P}}}
\widetilde{\omega}_{\mathbf{P}/\mathbf{A}}^{p-r})
\longrightarrow
\varGamma(\mathbf{A}_{\sigma}\cap\mathbf{P}_{\tau},
\widetilde{\omega}_{\mathbf{P}_{\tau}/\mathbf{A}}^{p-r}).
$$
We denote by
$$
\Phi_{\tau}:\omega_{\mathbf{P}/\mathbf{A}}^{r}(\log{D})
\otimes_{\mathcal{O}_{\mathbf{P}}}
\widetilde{\omega}_{\mathbf{P}/\mathbf{A}}^{p-r}\longrightarrow
\iota_{\tau*}\widetilde{\omega}_{\mathbf{P}_{\tau}/\mathbf{A}}^{p-r}
$$
the homomorphism of sheaves defined by $\Phi_{\sigma,\tau}$ for
$\sigma\in\Sigma$, and define a homomorphism by
$$
\Phi=\bigoplus_{\tau\in\Sigma_{\pi}(r)}\Phi_{\tau}:
\omega_{\mathbf{P}/\mathbf{A}}^{r}(\log{D})
\otimes_{\mathcal{O}_{\mathbf{P}}}
\widetilde{\omega}_{\mathbf{P}/\mathbf{A}}^{p-r}\longrightarrow
\bigoplus_{\tau\in\Sigma_{\pi}(r)}
\iota_{\tau*}\widetilde{\omega}_{\mathbf{P}_{\tau}/\mathbf{A}}^{p-r}.
$$
\begin{lemma}\label{poi}
 The kernel of the natural homomorphism
 $$
 \omega_{\mathbf{P}/\mathbf{A}}^{r}(\log{D})
 \otimes_{\mathcal{O}_{\mathbf{P}}}
 \widetilde{\omega}_{\mathbf{P}/\mathbf{A}}^{p-r}
 \longrightarrow\omega_{\mathbf{P}/\mathbf{A}}^{p}(\log{D});\
 \eta\otimes w\longmapsto\eta\wedge w
 $$
 is contained in the kernel of $\Phi$.
\end{lemma}
\begin{proof}
 Let $\sigma$ be a cone in $\Sigma$, let $u$ be an element in
 $M\cap\sigma^{\vee}$, and let
 $(w_{1},\ldots,w_{n})$ be a $k$-basis of $M_{\pi,k}$ such that
 $(w_{1},\ldots,w_{s})$ is a $k$-basis of $H_{\sigma,u}$.
 For $I\subset\{1,\ldots,n\}$ with $|I|=r$, and for
 $J\subset\{1,\ldots,s\}$ with $|J|=p-r$, vectors
 $w_{I}\in\bigwedge^{r}M_{\pi,k}$ and
 $w_{J}\in\bigwedge^{p-r}H_{\sigma,u}$ are defined by the same way as
 $(\ref{mi})$.
 If
 $$
 w=\sum_{|I|=r}\sum_{|J|=p-r}a_{I,J}w_{I}\otimes
 w_{J}\hspace*{20pt}(a_{I,J}\in k)
 $$
 is contained in the kernel of the natural homomorphism
 $$
 \bigwedge^{r}M_{\pi,k}\otimes_{k}\bigwedge^{p-r}H_{\sigma,u}
 \longrightarrow\bigwedge^{p}M_{\pi,k},
 $$
 then for $L\subset\{1,\ldots,n\}$ with $|L|=p$,
 we have
 \begin{equation}\label{wa}
  \sum_{I\cup J=L}\sgn{(I,J)}a_{I,J}=0,
 \end{equation}
 where $\sgn{(I,J)}$ is defined by
 $w_{I}\wedge w_{J}=\sgn{(I,J)}w_{L}$.
 For $\tau\in\Sigma_{\pi}(r)$ with $\tau\subset\sigma$, we have to prove
 $\Phi_{\sigma,\tau}(w\chi^{u})=0$.
 In the case $u\notin\tau^{\perp}$, it is clear, so we assume that
 $u\in\tau^{\perp}$.
 Let $I\subset\{1,\ldots,n\}$ be a subset satisfying $|I|=r$ and
 $\phi_{\tau}(w_{I})\neq0$.
 Then we have $I\cap\{1,\ldots,s\}=\emptyset$, because
 $\langle w_{i},1\otimes v\rangle=0$ for $1\leq i\leq s$ and
 $v\in N\cap\tau_{\mathbf{R}}$.
 By $(\ref{wa})$, we have $a_{I,J}=0$ for any $J\subset\{1,\ldots,s\}$,
 hence
 $$
 \Phi_{\sigma,\tau}(w\chi^{u})=
 \sum_{|J|=p-r}\sum_{|I|=r}
 a_{I,J}\phi_{\tau}(w_{I})w_{J}\chi^{u}_{\tau}=0.
 $$
\end{proof}
We denote by
$W_{r}\omega_{\mathbf{P}/\mathbf{A}}^{p}(\log{D})$
the image of the natural homomorphism
$\omega_{\mathbf{P}/\mathbf{A}}^{r}(\log{D})
\otimes_{\mathcal{O}_{\mathbf{P}}}
\widetilde{\omega}_{\mathbf{P}/\mathbf{A}}^{p-r}
\rightarrow\omega_{\mathbf{P}/\mathbf{A}}^{p}(\log{D})$.
By Lemma~$\ref{poi}$, the homomorphism $\Phi$ induces a homomorphism
$$
\mathrm{Res}:
W_{r}\omega_{\mathbf{P}/\mathbf{A}}^{p}(\log{D})
\longrightarrow\bigoplus_{\tau\in\Sigma_{\pi}(r)}
\iota_{\tau*}\widetilde{\omega}_{\mathbf{P}_{\tau}/\mathbf{A}}^{p-r},
$$
which is called Poincar\'{e} residue map.
The Poincar\'{e} residue map has the following fundamental property like
the case on $\varOmega_{\mathbf{P}}^{p}$ by \cite{d1} or \cite{da}.
\begin{theorem}\label{res}
 Let $\mathbf{A}$ be an affine toric variety, let
 $\mathbf{P}=\mathbf{P}_{\Sigma}$ be a simplicial toric variety such
 that $m(\Sigma)$ is prime to the characteristic of $k$, and let
 $\pi:\mathbf{P}\rightarrow\mathbf{A}$
 be an equivariant morphism.
 There is an exact sequence of $\mathcal{O}_{\mathbf{P}}$-modules
 $$
 0\longrightarrow
 W_{r-1}\omega_{\mathbf{P}/\mathbf{A}}^{p}(\log{D})
 \longrightarrow
 W_{r}\omega_{\mathbf{P}/\mathbf{A}}^{p}(\log{D})
 \overset{\mathrm{Res}}{\longrightarrow}
 \bigoplus_{\tau\in\Sigma_{\pi}(r)}
 \iota_{\tau*}\widetilde{\omega}_{\mathbf{P}_{\tau}/\mathbf{A}}^{p-r}
 \longrightarrow0.
 $$
\end{theorem}
\begin{proof}
 We check this on the affine coordinate $\mathbf{A}_{\sigma}$ for
 $\sigma\in\Sigma$.
 First we prove that the Poincar\'{e} residue map is surjective.
 Let $\tau\subset\sigma$ be a cone in $\Sigma_{\pi}(r)$ and let
 $v_{1},\ldots,v_{c}\in N$ be satisfying
 \begin{align*}
  \{v_{1},\ldots,v_{c}\}&=\{v_{\rho}\mid\rho\in\Sigma(1)\
  \text{and}\ \rho\subset\sigma\},\\
  \{v_{1},\ldots,v_{l}\}&=\{v_{\rho}\mid\rho\in\Sigma_{\pi}(1)\
  \text{and}\ \rho\subset\sigma\},\\
  \{v_{1},\ldots,v_{r}\}&=\{v_{\rho}\mid\rho\in\Sigma_{\pi}(1)\
  \text{and}\ \rho\subset\tau\}.
 \end{align*}
 Since $m(\sigma)$ is prime to the characteristic of $k$, the vectors
 $1\otimes v_{1},\ldots,1\otimes v_{l}$ are linearly independent in
 $N_{\pi,k}$, where $N_{\pi,k}$ denotes the kernel of
 $\pi_{*k}:k\otimes_{\mathbf{Z}}N\rightarrow k\otimes_{\mathbf{Z}}N'$.
 Then there exists a $k$-basis $(w_{1},\ldots,w_{n})$ of $M_{\pi,k}$
 such that
 $$
 \langle w_{i},1\otimes v_{j}\rangle=
 \begin{cases}
  1&(i=j),\\
  0&(i\neq j)
 \end{cases}
 $$
 for $1\leq i\leq n$ and $1\leq j\leq l$.
 For a subset $I\subset\{1,\ldots,l\}$ with $|I|=r$, we define a cone by
 $\tau_{I}=\sum_{i\in I}\mathbf{R}_{\geq0}v_{i}\in\Sigma_{\pi}(r)$.
 Then we have
 $$
 \phi_{\tau_{I}}(\eta)=
 \begin{cases}
  1&(I=\{1,\ldots,r\}),\\
  0&(I\neq\{1,\ldots,r\}),
 \end{cases}
 $$
 where $\eta=w_{1}\wedge\cdots\wedge w_{r}$.
 For $u\in M\cap\sigma\cap\tau^{\perp}$ and
 $w\in\bigwedge^{r-p}H_{\sigma,u}$, we have
 $w\chi^{u}_{\tau}=\Res{(\eta\wedge w\chi^{u})}$.\par
 Next we show that
 $W_{r-1}\omega_{\mathbf{P}/\mathbf{A}}^{p}(\log{D})$ is
 contained in the kernel of the Poincar\'{e} residue map.
 For $u\in M\cap\sigma^{\vee}$, we define a subset of $\{1,\ldots,l\}$
 by
 $$
 J_{u}=\{j\in\{1,\ldots,l\}\mid\langle u,v_{j}\rangle=0\}.
 $$
 Then $(w_{J};\ |J|=q,\ J\cap J_{u}=\emptyset)$ is a $k$-basis of
 $\bigwedge^{q}H_{\sigma,u}$.
 For $L\subset\{1,\ldots,n\}$ with $|L|=r-1$, and for
 $J\subset\{1,\ldots,n\}$ with $|J|=p-r+1$ and $J\cap J_{u}=\emptyset$,
 we have
 $$
 \Res{(w_{L}\wedge w_{J}\chi^{u})}=
 \begin{cases}
  \Phi(w_{L}\wedge w_{j_{1}}\otimes
  w_{J\smallsetminus\{j_{1}\}}\chi^{u})&(L\cap J=\emptyset),\\
  0&(L\cap J\neq\emptyset),
 \end{cases}
 $$
 where $j_{1}$ is the smallest integer in $J$.
 For $I\subset\{1,\ldots,l\}$ with $|I|=r$,
 if $u\in\tau_{I}^{\perp}$, then
 $\phi_{\tau_{I}}(w_{L}\wedge w_{j_{1}})=0$, so we have
 $\Res{(w_{L}\wedge w_{J}\chi^{u})}=0$.\par
 Finally we assume that
 $$
 \sum_{|I|=r}\sum_{|J|=p-r}a_{I,J}w_{I}\otimes
 w_{J}\chi^{u}\hspace*{20pt}(a_{I,J}\in k)
 $$
 is contained in the kernel of $\Phi$.
 For $L\subset\{1,\ldots,l\}$ with $|L|=r$,
 if $L\subset J_{u}$, then $u\in\tau_{L}^{\perp}$, hence
 $$
 a_{L,J}=\phi_{\tau_{L}}(\sum_{|I|=r}a_{I,J}w_{I})=0.
 $$
 We have
 \begin{align*}
  \sum_{|I|=r}\sum_{|J|=p-r}a_{I,J}w_{I}\wedge w_{J}
  &=\sum_{I\nsubseteq J_{u}}\sum_{|J|=p-r}a_{I,J}w_{I}\wedge w_{J}\\
  &=\sum_{I\nsubseteq J_{u}}
  \sgn{(I\smallsetminus\{i\},i)}
  \sum_{|J|=p-r}a_{I,J}w_{I\smallsetminus\{i\}}
  \wedge w_{i}\wedge w_{J},
 \end{align*}
 where we take $i\in I\smallsetminus J_{u}$ for each
 $I\nsubseteq J_{u}$.
 Since $(J\cup\{i\})\cap J_{u}=\emptyset$,
 this is contained in
 $\varGamma(\mathbf{A}_{\sigma},
 W_{r-1}\omega_{\mathbf{P}/\mathbf{A}}^{p}(\log{D}))$.
\end{proof}
\begin{proposition}\label{sj}
 If $\pi:\mathbf{P}\rightarrow\mathbf{A}$ is proper, then for any ample
 invertible sheaf $\mathcal{L}$ on $\mathbf{P}$, the global Poincar\'{e}
 residue map
 $$
 H^{0}(\mathbf{P},W_{r}\omega_{\mathbf{P}/S}^{p}(\log{D})
 \otimes_{\mathcal{O}_{\mathbf{P}}}\mathcal{L})
 \overset{\mathrm{Res}}{\longrightarrow}
 \bigoplus_{\tau\in\Sigma_{\pi}(r)}
 H^{0}(\mathbf{P}_{\tau},
 \widetilde{\omega}_{\mathbf{P}_{\tau}/\mathbf{A}}^{p-r}
 \otimes_{\mathcal{O}_{\mathbf{P}_{\tau}}}
 \iota_{\tau}^{*}\mathcal{L})
 $$
 is surjective.
\end{proposition}
\begin{proof}
 There is a $\Sigma$-linear support function $h$ such that
 $\mathcal{L}\simeq\mathcal{O}_{\mathbf{P}}(D_{h})$.
 We define a subset of $\Sigma_{\pi}(1)$ by
 $$
 \Sigma_{\pi,u=h}(1)=\{\rho\in\Sigma_{\pi}(1)\mid
 \langle u,v_{\rho}\rangle=h(v_{\rho})\},
 $$
 and define a $k$-subspace of $M_{\pi,k}$ by
 $$
 H_{u,h}=\{w\in M_{\pi,\mathbf{R}}\mid
 \langle w,1\otimes v_{\rho}\rangle=0\
 \text{for any}\ \rho\in\Sigma_{\pi,u=h}(1)\}.
 $$
 By Lemma~$\ref{li}$, there is a natural isomorphism
 $$
 H^{0}(\mathbf{P},\widetilde{\omega}_{\mathbf{P}/S}^{p-r}
 \otimes_{\mathcal{O}_{\mathbf{P}}}\mathcal{O}_{\mathbf{P}}(D_{h}))
 \simeq\bigoplus_{u\in M\cap\Delta_{h}}\bigwedge^{p-r}H_{u,h}\cdot\chi^{u},
 $$
 where
 $$
 \Delta_{h}=\Delta_{D_{h}}=
 \{u\in M_{\mathbf{R}}\mid\langle u,v\rangle\geq h(v)\
 \text{for any}\ v\in|\Sigma|\}.
 $$
 Let $\tau$ be a cone in $\Sigma_{\pi}(r)$.
 There is a linear function $u_{\tau}\in M$ which coincides with $h$ on
 $\tau$.
 Then $\Sigma$-linear support function $h-u_{\tau}$ induces a
 $\Sigma_{\tau}$-linear support function
 $h_{\tau}:N_{\tau,\mathbf{R}}\rightarrow\mathbf{R}$,
 and there is an
 isomorphism
 $\iota_{\tau}^{*}\mathcal{O}_{\mathbf{P}}(D_{h-u_{\tau}})\simeq
 \mathcal{O}_{\mathbf{P}_{\tau}}(D_{h_{\tau}})$.
 Since $\pi$ is proper and $h$ is upper convex, there is a natural
 identification
 $$
 M\cap\Delta_{h-u_{\tau}}\cap\tau^{\perp}\simeq
 M_{\tau}\cap\Delta_{h_{\tau}}.
 $$
 Since $m(\Sigma)$ is prime to the characteristic of $k$, for
 $u\in M\cap\Delta_{h-u_{\tau}}\cap\tau^{\perp}$, there is a natural
 identification
 $$
 H_{u,h-u_{\tau}}\simeq
 \{w\in(M_{\tau})_{\pi_{\tau},u}\mid
 \langle w,1\otimes v_{\rho}\rangle=0\ \text{for any}\
 \rho\in(\Sigma_{\tau})_{\pi_{\tau},u=h_{\tau}}(1)\}.
 $$
 By the isomorphism
 $$
 \mathcal{O}_{\mathbf{P}}(D_{h})\simeq
 \mathcal{O}_{\mathbf{P}}(D_{h-u_{\tau}});\
 \chi^{u}\longleftrightarrow\chi^{u-u_{\tau}},
 $$
 the homomorphism
 $$
 \Phi_{\tau}:
 \varGamma(\mathbf{P},\omega_{\mathbf{P}/\mathbf{A}}^{r}(\log{D})
 \otimes\widetilde{\omega}_{\mathbf{P}/\mathbf{A}}^{p-r}
 \otimes\mathcal{O}_{\mathbf{P}}(D_{h}))
 \longrightarrow
 \varGamma(\mathbf{P}_{\tau},
 \widetilde{\omega}_{\mathbf{P}_{\tau}/\mathbf{A}}^{p-r}
 \otimes\mathcal{O}_{\mathbf{P}_{\tau}}(D_{h_{\tau}}))
 $$
 is given by
 \begin{align*}
  \Phi_{\tau}:
  \bigwedge^{r}M_{\pi,k}&\otimes_{k}
  \bigoplus_{u\in M\cap\Delta_{h}}\bigwedge^{p-r}H_{u,h}\cdot\chi^{u}
  &\longrightarrow&
  \bigoplus_{u\in M\cap\Delta_{h-u_{\tau}}\cap\tau^{\perp}}
  \bigwedge^{p-r}H_{u,h-u_{\tau}}\cdot\chi^{u}_{\tau};\\
  &\eta\otimes w\chi^{u}
  &\longmapsto&
  \begin{cases}
   \phi_{\tau}(\eta)w\chi^{u-u_{\tau}}_{\tau}
   &(u-u_{\tau}\in\tau^{\perp}),\\
   0&(u-u_{\tau}\notin\tau^{\perp}),
  \end{cases}
 \end{align*}
 where we remark that $H_{u,h}=H_{u-u_{\tau},h-u_{\tau}}$.
 Since $h$ is strictly upper convex, for
 $u\in M\cap\Delta_{h-u_{\tau}}\cap\tau^{\perp}$, there exists a cone
 $\sigma\in\Sigma$ such that $\sigma(1)\supset\Sigma_{\pi,u=h}(1)$.
 Using this cone $\sigma$, we give a $k$-basis $(w_{1},\ldots,w_{n})$ of
 $M_{\pi,\tau}$ in the same way as the proof of Theorem~$\ref{res}$,
 and define a vector in $\bigwedge^{r}M_{\pi,k}$ by
 $\eta=w_{1}\wedge\cdots\wedge w_{r}$.
 Then for $w\in \bigwedge^{p-r}H_{u,h-u_{\tau}}$, we have
 $w\chi^{u}=\Phi(\eta\otimes w\chi^{u+u_{\tau}})$.
\end{proof}
We prove a vanishing theorem of Bott type for the cohomology of
the sheaf of relative logarithmic differential forms, that is reduced to
the following vanishing theorem for invertible sheaves on toric
varieties.
The idea is same as \cite{bc}.
\begin{theorem}[\cite{f}]\label{vl}
 Let $\mathbf{P}_{\Sigma}$ be a toric variety, and let $\mathcal{L}$ be
 an invertible sheaf on $\mathbf{P}_{\Sigma}$.
 If the support $|\Sigma|$ is convex, and $\mathcal{L}$ is generated by
 global sections, then for $q\geq1$,
 $$
 H^{q}(\mathbf{P}_{\Sigma},\mathcal{L})=0.
 $$
\end{theorem}
\begin{theorem}\label{van}
 Let $\mathbf{A}$ be an affine toric variety, let $\mathbf{P}$ be a
 simplicial toric variety such that $m(\Sigma)$ is prime to the
 characteristic of $k$, and let
 $\pi:\mathbf{P}\rightarrow\mathbf{A}$ be a proper equivariant morphism,
 and let $\mathcal{L}$ be an invertible sheaf on $\mathbf{P}$.
 \renewcommand{\labelenumi}{$(\arabic{enumi})$\hspace*{5pt}}
 \begin{enumerate}
  \item If $\mathcal{L}$ is generated by global sections,
	then for $p\geq r\geq0$ and $q\geq p-r+1$,
	$$
	H^{q}(\mathbf{P},W_{r}\omega_{\mathbf{P}/\mathbf{A}}^{p}
	(\log{D})\otimes_{\mathcal{O}_{\mathbf{P}}}\mathcal{L})=0.
	$$
  \item If $\mathcal{L}$ is ample, then for $p\geq r\geq0$ and $q\geq1$,
	$$
	H^{q}(\mathbf{P},W_{r}\omega_{\mathbf{P}/\mathbf{A}}^{p}
	(\log{D})\otimes_{\mathcal{O}_{\mathbf{P}}}\mathcal{L})=0.
	$$
 \end{enumerate}
\end{theorem}
\begin{proof}
 We prove this by induction on $p-r$.
 If $p-r=0$, then
 $W_{p}\omega_{\mathbf{P}/\mathbf{A}}^{p}(\log{D})$
 is isomorphic to the free $\mathcal{O}_{\mathbf{P}}$-module
 $\mathcal{O}_{\mathbf{P}}\otimes_{k}\bigwedge^{p}M_{\pi,k}$.
 By Theorem~$\ref{vl}$, we have
 $H^{q}(W_{p}\omega_{\mathbf{P}/\mathbf{A}}^{p}(\log{D})
 \otimes\mathcal{L})=0$.\par
 For an integer $l>0$, we assume that Theorem~$\ref{van}$ is true for
 $p-r<l$.
 Then we prove Theorem~$\ref{van}$ for $p\geq r\geq0$ with $p-r=l$.
 By Theorem~$\ref{res}$, there is an exact sequence
\begin{multline*}
 H^{q-1}(W_{r+1}\omega_{\mathbf{P}/\mathbf{A}}^{p}(\log{D})
 \otimes\mathcal{L})
 \longrightarrow
 \bigoplus_{\tau\in\Sigma_{\pi}(r+1)}
 H^{q-1}(\widetilde{\omega}_{\mathbf{P}_{\tau}/\mathbf{A}}^{p-r-1}
 \otimes\iota_{\tau}^{*}\mathcal{L})\\
 \longrightarrow
 H^{q}(W_{r}\omega_{\mathbf{P}/\mathbf{A}}^{p}(\log{D})
 \otimes\mathcal{L})
 \longrightarrow
 H^{q}(W_{r+1}\omega_{\mathbf{P}/\mathbf{A}}^{p}(\log{D})
 \otimes\mathcal{L}).
\end{multline*}
 By the assumption for induction, we have
 $H^{q}(W_{r+1}\omega_{\mathbf{P}/\mathbf{A}}^{p}(\log{D})
 \otimes\mathcal{L})=0$
 for $q\geq p-r$, and
 $H^{q-1}(\widetilde{\omega}_{\mathbf{P}_{\tau}/\mathbf{A}}^{p-r-1}
 \otimes\iota_{\tau}^{*}\mathcal{L})=0$
 for $q\geq p-r+1$.
 Hence we have
 $H^{q}(W_{r}\omega_{\mathbf{P}/\mathbf{A}}^{p}(\log{D})
 \otimes\mathcal{L})=0$
 for $q\geq p-r+1$.\par
 In the case when $\mathcal{L}$ is ample, by the assumption for
 induction, we have
 $H^{q}(W_{r+1}\omega_{\mathbf{P}/\mathbf{A}}^{p}(\log{D})
 \otimes\mathcal{L})=0$
 for $q\geq1$, and
 $H^{q-1}(\widetilde{\omega}_{\mathbf{P}_{\tau}/\mathbf{A}}^{p-r-1}
 \otimes\iota_{\tau}^{*}\mathcal{L})=0$
 for $q\geq 2$.
 By Proposition~$\ref{sj}$, the homomorphism
 $$
 H^{0}(W_{r+1}\omega_{\mathbf{P}/\mathbf{A}}^{p}(\log{D})
 \otimes\mathcal{L})
 \longrightarrow
 \bigoplus_{\tau\in\Sigma_{\pi}(r+1)}
 H^{0}(\widetilde{\omega}_{\mathbf{P}_{\tau}/\mathbf{A}}^{p-r-1}
 \otimes\iota_{\tau}^{*}\mathcal{L})
 $$
 is surjective.
 Hence we have
 $H^{q}(W_{r}\omega_{\mathbf{P}/\mathbf{A}}^{p}(\log{D})
 \otimes\mathcal{L})=0$
 for $q\geq1$.
\end{proof}
\begin{corollary}\label{v}
 \renewcommand{\labelenumi}{$(\arabic{enumi})$\hspace*{5pt}}
 \begin{enumerate}
  \item If $\mathcal{L}$ is generated by global sections,
	then for $q\geq p+1$,
	$$
	H^{q}(\mathbf{P},
	\widetilde{\omega}_{\mathbf{P}/\mathbf{A}}^{p}
	\otimes_{\mathcal{O}_{\mathbf{P}}}\mathcal{L})=0.
	$$
  \item If $\mathcal{L}$ is ample, then for $q\geq1$,
	$$
	H^{q}(\mathbf{P},
	\widetilde{\omega}_{\mathbf{P}/\mathbf{A}}^{p}
	\otimes_{\mathcal{O}_{\mathbf{P}}}\mathcal{L})=0.
	$$
 \end{enumerate}
\end{corollary}
In the rest of this section, we prove the Euler exact sequence for the
sheaf of relative logarithmic differential forms.
Here we introduce the notion of logarithmically smoothness.
\begin{definition}\label{logsm}
 If the rank of the locally free
 $\mathcal{O}_{\mathbf{P}^{\mathrm{reg}}}$-module
 $\omega_{\mathbf{P}^{\mathrm{reg}}/\mathbf{A}}^{1}$
 is equal to
 $\dim{\mathbf{P}}-\dim{\mathbf{A}}$, then
 we call that $\pi$ is {\itshape logarithmically smooth}.
\end{definition}
\begin{remark}
 The following conditions are equivalent:
 \renewcommand{\labelenumi}{$(\arabic{enumi})$\hspace*{5pt}}
 \begin{enumerate}
  \item $\pi$ is logarithmically smooth in the sense of
	Definition~$\ref{logsm}$;
  \item $\pi\circ j:(\mathbf{P}^{\mathrm{reg}},
	\mathcal{N}_{\mathbf{P}^{\mathrm{reg}}}^{E})
	\rightarrow
	(\mathbf{A},\mathcal{N}_{\mathbf{A}}^{\mathrm{can}})$
	is logarithmically smooth in the sense of \cite{k};
  \item $\pi:(\mathbf{P},
	\mathcal{N}_{\mathbf{P}}^{\mathrm{can}})
	\rightarrow
	(\mathbf{A},\mathcal{N}_{\mathbf{A}}^{\mathrm{can}})$
	is logarithmically smooth in the sense of \cite{k};
  \item the cokernel of $\pi_{*}:N\rightarrow N'$ is finite, whose order
	is prime to the characteristic of $k$.
  \item $\pi^{*}:M'\rightarrow M$ is injective and the order of the
	torsion part of the cokernel of $\pi^{*}$ is prime to the
	characteristic of $k$.
 \end{enumerate}
\end{remark}
We denote by $\widehat{N}_{\pi}$ the free $\mathbf{Z}$-module
generated by $\Sigma_{\pi}(1)$, and denote by $P_{\pi}$ the kernel of
the homomorphism
$$
\psi:\widehat{N}_{\pi}\longrightarrow N;\
\sum_{\rho\in\Sigma_{\pi}(1)}b_{\rho}\rho\longmapsto
\sum_{\rho\in\Sigma_{\pi}(1)}b_{\rho}v_{\rho}.
$$
We identify the group of $\mathbf{T}$-invariant Weil divisors on
$\mathbf{P}\smallsetminus E$ with the dual $\mathbf{Z}$-module
$\widehat{M}_{\pi}$ of $\widehat{N}_{\pi}$ by the pairing
$$
\langle\sum_{\rho\in\Sigma_{\pi}(1)}a_{\rho}\mathbf{P}_{\rho},
\sum_{\rho\in\Sigma_{\pi}(1)}b_{\rho}\rho\rangle=
\sum_{\rho\in\Sigma_{\pi}(1)}a_{\rho}b_{\rho}.
$$
Then the divisor class group $\Cl{(\mathbf{P}\smallsetminus E)}$ is
naturally isomorphic to the cokernel of the dual homomorphism
$\psi^{*}:M\longrightarrow\widehat{M}_{\pi}$.
\begin{theorem}\label{eu}
 Let $\mathbf{A}$ be an affine toric variety, let $\mathbf{P}$ be a
 simplicial toric variety such that $m(\Sigma)$ is prime to the
 characteristic of $k$, and let
 $\pi:\mathbf{P}\rightarrow\mathbf{A}$ be a logarithmically smooth
 proper equivariant morphism.
 Then there is an exact sequence of $\mathcal{O}_{\mathbf{P}}$-modules
 $$
 0\longrightarrow
 \widetilde{\omega}_{\mathbf{P}/\mathbf{A}}^{1}
 \overset{\Psi^{*}}{\longrightarrow}
 \bigoplus_{\rho\in\Sigma_{\pi}(1)}
 \mathcal{O}_{\mathbf{P}}(-\mathbf{P}_{\rho})
 \overset{\gamma}{\longrightarrow}
 \mathcal{O}_{\mathbf{P}}\otimes_{\mathbf{Z}}
 \Cl{(\mathbf{P}\smallsetminus E)}
 \longrightarrow0.
 $$
\end{theorem}
\begin{proof}
 We denote by $N_{\pi}$ the kernel of $\pi_{*}:N\rightarrow N'$.
 Since $\pi$ is proper, and $m(\Sigma)$ is prime to the characteristic
 of $k$, the cokernel of $\psi:\widehat{N}_{\pi}\rightarrow N_{\pi}$
 is finite, whose order is prime to the characteristic of $k$.
 So we have an exact sequence of $k$-vector spaces
 $$
 0\longrightarrow k\otimes_{\mathbf{Z}}P_{\pi}
 \longrightarrow k\otimes_{\mathbf{Z}}\widehat{N}_{\pi}
 \overset{\psi_{k}}{\longrightarrow}k\otimes_{\mathbf{Z}}N_{\pi}
 \longrightarrow0.
 $$
 Since $\pi$ is logarithmically smooth, we have
 $M_{\pi,k}\simeq
 \Hom_{k}
 {(k\otimes_{\mathbf{Z}} N_{\pi},k)}$.
 Hence there is an exact sequence of $k$-vector spaces
 $$
 0\longrightarrow M_{\pi,k}
 \overset{\psi_{k}^{*}}{\longrightarrow}
 k\otimes_{\mathbf{Z}}\widehat{M}_{\pi}
 \longrightarrow k\otimes_{\mathbf{Z}}
 \Hom_{\mathbf{Z}}{(P_{\pi},\mathbf{Z})}
 \longrightarrow0,
 $$
 which induce an exact sequence of $\mathcal{O}_{\mathbf{P}}$-modules
 $$
 0\longrightarrow
 \omega_{\mathbf{P}/\mathbf{A}}^{1}(\log{D})
 \overset{\psi_{\mathcal{O}_{\mathbf{P}}}^{*}}{\longrightarrow}
 \mathcal{O}_{\mathbf{P}}\otimes_{\mathbf{Z}}\widehat{M}_{\pi}
 \longrightarrow
 \mathcal{O}_{\mathbf{P}}\otimes_{\mathbf{Z}}
 \Cl{(\mathbf{P}\smallsetminus E)}
 \longrightarrow0.
 $$
 The exact sequence in Theorem~$\ref{eu}$ is shown by the commutative
 diagram
 $$
 \begin{CD}
  &&0&&0&&&&\\
  &&\downarrow&&\downarrow&&&&\\
  &&
  \widetilde{\omega}_{\mathbf{P}/\mathbf{A}}^{1}
  &&
  \displaystyle{\bigoplus_{\rho\in\Sigma_{\pi}(1)}}
  \mathcal{O}_{\mathbf{P}}(-\mathbf{P}_{\rho})
  &&&&\\
  &&\downarrow&&\downarrow&&&&\\
  0&\longrightarrow&
  \omega_{\mathbf{P}/\mathbf{A}}^{1}(\log{D})
  &\overset{\psi_{\mathcal{O}_{\mathbf{P}}}^{*}}{\longrightarrow}&
  \mathcal{O}_{\mathbf{P}}\otimes_{\mathbf{Z}}\widehat{M}_{\pi}
  &\longrightarrow&
  \mathcal{O}_{\mathbf{P}}\otimes_{\mathbf{Z}}
  \Cl{(\mathbf{P}\smallsetminus E)}
  &\longrightarrow&0\\
  &&\mathrm{Res}\downarrow\hspace*{23pt}&&\downarrow&&&&\\
  &&
  \displaystyle{\bigoplus_{\rho\in\Sigma_{\pi}(1)}}
  \iota_{\rho*}\mathcal{O}_{\mathbf{P}_{\rho}}
  &=&
  \displaystyle{\bigoplus_{\rho\in\Sigma_{\pi}(1)}}
  \iota_{\rho*}\mathcal{O}_{\mathbf{P}_{\rho}}&&&&\\
  &&\downarrow&&\downarrow&&&&\\
  &&0&&0,&&&&
 \end{CD}
 $$
 where the left vertical sequence is proved in Theorem~$\ref{res}$.
\end{proof}
\section{Hypersurfaces in toric varieties}\label{s4}
Let $\mathbf{A}=\Spec{A}$ be an affine toric variety with a torus
invariant point $0$, let $\mathbf{P}=\mathbf{P}_{\Sigma}$ be a
nonsingular toric variety, and let
$\pi:\mathbf{P}\rightarrow\mathbf{A}$ be a logarithmically smooth
proper equivariant morphism.
Then we remark that $\dim{\mathbf{P}}=\dim{|\Sigma|}$, and
$\omega_{\mathbf{P}/\mathbf{A}}^{1}
=\omega_{\mathbf{P}^{\mathrm{reg}}/\mathbf{A}}^{1}
=\widetilde{\omega}_{\mathbf{P}/\mathbf{A}}^{1}$
is a locally free $\mathcal{O}_{\mathbf{P}}$-module of rank
$n=\dim{\mathbf{P}}-\dim{\mathbf{A}}$.\par
The homogeneous coordinate ring of
$\mathbf{P}$ is defined in \cite{c} as a $\Cl{(\mathbf{P})}$-graded
polynomial ring
$$
S_{\mathbf{P}}=k[z_{\rho};\ \rho\in\Sigma(1)]
=\bigoplus_{\beta\in\Cl{(\mathbf{P})}}S_{\mathbf{P}}^{\beta}
$$
with
$$
\deg{z_{\rho}}=[\mathbf{P}_{\rho}]\in\Cl{(\mathbf{P})}.
$$
For a $\mathbf{T}_{N}$-invariant divisor
$B=\sum_{\rho\in\Sigma(1)}b_{\rho}\mathbf{P}_{\rho}$, there is a
natural isomorphism
\begin{equation}\label{hm}
 S_{\mathbf{P}}^{[B]}
  \simeq
  H^{0}(\mathbf{P},\mathcal{O}_{\mathbf{P}}(B));\
  \prod_{\rho\in\Sigma(1)}z_{\rho}^{\langle u,v_{\rho}\rangle+b_{\rho}}
  \leftrightarrow\chi^{u}.
\end{equation}
By the $A$-module structure on
$H^{0}(\mathbf{P},\mathcal{O}_{\mathbf{P}}(B))
=H^{0}(\mathbf{A},\pi_{*}\mathcal{O}_{\mathbf{P}}(B))$,
the homogeneous coordinate ring $S_{\mathbf{P}}$ has an $A$-algebra
structure.\par
Let $X$ be a hypersurface in $\mathbf{P}$.
Then there is a $\mathbf{T}_{N}$-invariant divisor
$B=\sum_{\rho\in\Sigma(1)}b_{\rho}\mathbf{P}_{\rho}$ such that
$$
[X]=[B]\in\Cl{(\mathbf{P})}.
$$
Using the isomorphism $(\ref{hm})$, the hypersurface $X$ is defined by a
$\Cl{(\mathbf{P})}$-homogeneous polynomial $F$, which does not depend on
the choice of $B$ modulo $k^{\times}$.
We define the Jacobian ring of $X$ over $\mathbf{A}$ by
$$
R_{X/\mathbf{A}}=
S_{\mathbf{P}}/
(\frac{\partial F}{\partial z_{\rho}};\ \rho\in\Sigma_{\pi}(1)),
$$
which is a $\Cl{(\mathbf{P})}$-graded $A$-algebra.
\begin{remark}
 If $\mathbf{A}=\Spec{k}$, then $\mathbf{P}$ is a complete toric
 variety, and $\Sigma_{\pi}(1)=\Sigma(1)$.
 In this case, our definition of Jacobian ring is same as \cite{bc}.
\end{remark}
We denote by $\mathcal{I}_{X/\mathbf{P}}$ the ideal sheaf of $X$ in
$\mathbf{P}$, and define a coherent $\mathcal{O}_{X}$-module by
$$
\omega_{X/\mathbf{A}}^{p}
=\Coker{(\mathcal{I}_{X/\mathbf{P}}/\mathcal{I}_{X/\mathbf{P}}^{2}
\otimes\omega_{\mathbf{P}/\mathbf{A}}^{p-1}\vert_{X}
\longrightarrow\omega_{\mathbf{P}/\mathbf{A}}^{p}\vert_{X};\
[f]\otimes\omega\mapsto (df\vert_{X}\wedge\omega))},
$$
where
$
d:\mathcal{O}_{\mathbf{P}}\rightarrow
\omega_{\mathbf{P}/\mathbf{A}}^{1}
$
is the differential operator given by
$$
d:
\bigoplus_{u\in M\cap\sigma^{\vee}}k\cdot\chi^{u}\longrightarrow
\bigoplus_{u\in M\cap\sigma^{\vee}}H_{\sigma,u}\cdot\chi^{u};\
\chi^{u}\longmapsto[1\otimes u]\chi^{u}.
$$
\begin{definition}
 Let $U$ be an open subset of $\mathbf{A}$.
 If $j_{U}^{*}\omega_{X/\mathbf{A}}^{1}$ is a locally free
 $\mathcal{O}_{X_{U}}$-module of rank $n-1$, then we call that $X$ is
 {\itshape logarithmic smooth} over $U$, where
 $j_{U}:X_{U}\rightarrow X$ denotes the open immersion
 $U\times_{\mathbf{A}}X\rightarrow X$.
\end{definition}
We define a coherent $\mathcal{O}_{\mathbf{P}}$-module by
$$
\omega_{\mathbf{P}/\mathbf{A}}^{p}(\log{X})
=\Ker{(\omega_{\mathbf{P}/\mathbf{A}}^{p}(X)
\longrightarrow\omega_{X/\mathbf{A}}^{p}\otimes
\mathcal{O}_{\mathbf{P}}(X)\vert_{X})}.
$$
If $X$ is ample in $\mathbf{P}$, then by the vanishing theorem
(Corollary~$\ref{v}$), we can calculate the cohomology of the sheaf
$\omega_{\mathbf{P}/\mathbf{A}}^{p}(\log{X})$,
using the following resolution.
\begin{lemma}\label{resol}
 If $X$ is logarithmically smooth over $U$, then for $0\leq p\leq n-2$,
 the following sequence is exact on
 $\mathbf{P}_{U}=U\otimes_{\mathbf{A}}{\mathbf{P}}$;
 \begin{multline*}
  0\longrightarrow
  \omega_{\mathbf{P}/\mathbf{A}}^{p+1}(\log{X})
  \longrightarrow
  \omega_{\mathbf{P}/\mathbf{A}}^{p+1}(X)
  \longrightarrow
  \omega_{\mathbf{P}/\mathbf{A}}^{p+2}(2X)\vert_{X}
  \longrightarrow\cdots\\
  \cdots\longrightarrow
  \omega_{\mathbf{P}/\mathbf{A}}^{n-1}((n-p-1)X)\vert_{X}
  \longrightarrow
  \omega_{\mathbf{P}/\mathbf{A}}^{n}((n-p)X)\vert_{X}
  \longrightarrow0.
 \end{multline*}
\end{lemma}
\begin{proof}
 We have the long exact sequence, by connecting the short exact sequence
 $$
 0\longrightarrow
 \mathcal{I}_{X/\mathbf{P}}/\mathcal{I}_{X/\mathbf{P}}^{2}
 \otimes\omega_{X/\mathbf{A}}^{p-1}
 \longrightarrow\omega_{\mathbf{P}/\mathbf{A}}^{p}\vert_{X}
 \longrightarrow\omega_{X/\mathbf{A}}^{p}
 \longrightarrow0
 $$
 on $\mathbf{P}_{U}$.
\end{proof}
The following is the main theorem in this paper, which describe the
cohomology of the sheaf of relative logarithmic forms by using the
Jacobian ring.
\begin{theorem}\label{main}
 If $X$ is ample and logarithmic smooth over an affine open subvariety
 $U=\Spec{A_{U}}$ of $\mathbf{A}=\Spec{A}$, and the class
 $[X]\in\Cl{(\mathbf{P}\smallsetminus E)}$ is not divisible by
 the characteristic of $k$, then for $0\leq p\leq n-1$, there is a
 natural isomorphism of $A_{U}$-modules
 $$
 H^{n-p-1}(\mathbf{P}_{U},
 \omega_{\mathbf{P}/\mathbf{A}}^{p+1}(\log{X}))
 \simeq
 A_{U}\otimes_{A}R_{X/\mathbf{A}}^{[(n-p)X-D]}.
 $$
\end{theorem}
We show some lemma for the proof of Theorem~$\ref{main}$.
\begin{lemma}\label{c1}
 The following diagram is commutative;
 $$
 \begin{CD}
  \mathcal{O}_{\mathbf{P}}\otimes_{\mathbf{Z}}
  \Cl{(\mathbf{P}\smallsetminus E)}^{*}
  &\overset{[X]}{\longrightarrow}&\mathcal{O}_{\mathbf{P}}\\
  \gamma^{*}\downarrow\hspace*{15pt}&&\hspace*{15pt}\downarrow F\\
  \displaystyle{\bigoplus_{\rho\in\Sigma_{\pi}(1)}}
  \mathcal{O}_{\mathbf{P}}
  (\mathbf{P}_{\rho})&
  \underset
  {(\frac{\partial F}{\partial z_{\rho}})_{\rho\in\Sigma_{\pi}(1)}}
  {\longrightarrow}&\mathcal{O}_{\mathbf{P}}(X),
 \end{CD}
 $$
 where the map $\gamma$ is defined in {\rm Theorem}~$\ref{eu}$.
\end{lemma}
\begin{proof}
 We denote by $F=\sum_{e}a_{e}z^{e}$ the $\Cl{(\mathbf{P})}$-homogeneous
 polynomial defining $X$, where $z^{e}$ is the monomial
 $\prod_{\rho\in\Sigma(1)}z^{e_{\rho}}$ of degree
 $[X]\in\Cl{(\mathbf{P})}$.
 Let $\varphi:\Cl{(\mathbf{P}\smallsetminus E)}\rightarrow\mathbf{Z}$ be
 a homomorphism.
 Then the image of $1\otimes\varphi$ by the map
 $(\frac{\partial F}{\partial z_{\rho}})_{\rho\in\Sigma_{\pi}(1)}
 \circ\gamma^{*}$
 is
 $$
 \sum_{\rho\in\Sigma_{\pi}(1)}\varphi([\mathbf{P}_{\rho}])z_{\rho}
 \frac{\partial F}{\partial z_{\rho}}
 =\sum_{\rho\in\Sigma_{\pi}(1)}\varphi([\mathbf{P}_{\rho}])
 \sum_{e}a_{e}e_{\rho}z^{e}
 =\sum_{e}a_{e}\varphi
 ([\sum_{\rho\in\Sigma_{\pi}(1)}e_{\rho}\mathbf{P}_{\rho}])z^{e}.
 $$
 Since
 $[\sum_{\rho\in\Sigma_{\pi}(1)}e_{\rho}\mathbf{P}_{\rho}]=[X]\in
 \Cl{(\mathbf{P}\smallsetminus E)}$,
 this is equal to
 $$
 \sum_{e}a_{e}\varphi([X])z^{e}=
 \varphi([X])F,
 $$
 which is the image of $1\otimes\varphi$ by the map
 $F\circ[X]$.
\end{proof}
\begin{lemma}\label{c2}
 The following diagram is commutative;
 $$
 \begin{CD}
  \displaystyle{\bigoplus_{\rho\in\Sigma_{\pi}(1)}}
  \omega_{\mathbf{P}/\mathbf{A}}^{n}
  (\mathbf{P}_{\rho})&
  \overset
  {(\frac{\partial F}{\partial z_{\rho}})_{\rho\in\Sigma_{\pi}(1)}}
  {\longrightarrow}
  &\omega_{\mathbf{P}/\mathbf{A}}^{n}(X)\\
  \Psi\downarrow\hspace*{15pt}&&\downarrow\\
  \omega_{\mathbf{P}/\mathbf{A}}^{n-1}
  &\underset{\delta}
  {\longrightarrow}&
  \omega_{\mathbf{P}/\mathbf{A}}^{n}(X)\vert_{X},
 \end{CD}
 $$
 where the map $\Psi$ is defined in {\rm Theorem}~$\ref{eu}$, and
 $\delta$ is defined as the restriction
 $$
 \omega_{\mathbf{P}/\mathbf{A}}^{n-1}\longrightarrow
 \omega_{X/\mathbf{A}}^{n-1}\simeq
 \omega_{\mathbf{P}/\mathbf{A}}^{n}(X)\vert_{X}.
 $$
\end{lemma}
\begin{proof}
 We check this on a local affine coordinate $\mathbf{A}_{\sigma}$ for
 $\sigma\in\Sigma$.
 There is a $\mathbf{T}$-invariant divisor
 $B=\sum_{\rho\in\Sigma(1)}b_{\rho}\mathbf{P}_{\rho}$
 such that $[B]=[X]\in\Cl{(\mathbf{P})}$, and $b_{\rho}=0$ for
 $\rho\nsubseteq\sigma$.
 Let $w$ be a $k[M\cap\sigma^{\vee}]$-basis of
 $\varGamma(\mathbf{A}_{\sigma},\omega_{\mathbf{P}/\mathbf{A}}^{n})$,
 and let $\chi^{u_{\rho}}$ be a $k[M\cap\sigma^{\vee}]$-basis of
 $\varGamma(\mathbf{A}_{\sigma},
 \mathcal{O}_{\mathbf{P}}(-\mathbf{P}_{\rho}))$.
 The image of
 $w\chi^{-u_{\rho}}\in
 \varGamma(\mathbf{A},
 \omega_{\mathbf{P}/\mathbf{A}}^{n}(\mathbf{P}_{\rho}))$
 by the map $\delta\circ\Psi$ is
 $$
 \frac{w}{f}\sum_{u\in M\cap\sigma^{\vee}}a_{u}
 \langle u,v_{\rho}\rangle\chi^{u-u_{\rho}}\vert_{X},
 $$
 where $f=\sum_{u\in M\cap\sigma^{\vee}}a_{u}\chi^{u}$
 is the local equation of $X$.
 On the other hand, the image of
 $w\chi^{-u_{\rho}}\in
 \varGamma(\mathbf{A},
 \omega_{\mathbf{P}/\mathbf{A}}^{n}(\mathbf{P}_{\rho}))$
 by the map $\frac{\partial F}{\partial z_{\rho}}$ is
 $$
 \frac{w}{f}\sum_{u}a_{u}(\langle u,v_{\rho}\rangle
 +b_{\rho})\chi^{u-u_{\rho}}=
 \frac{w}{f}(\sum_{u}a_{u}\langle u,v_{\rho}\rangle\chi^{u-u_{\rho}}
 +b_{\rho}f\chi^{-u_{\rho}}),
 $$
 whose restriction to $X$ is
 $(\delta\circ\Psi)(w\chi^{-u_{\rho}})\in
 \varGamma(\mathbf{A}_{\sigma},
 \omega_{\mathbf{P}/\mathbf{A}}^{n}(X)\vert_{X})$.
\end{proof}
\begin{proof}[Proof of Theorem~$\ref{main}$]
 In the case $p=n-1$, we have
 $$
 A_{U}\otimes_{A}R_{X/\mathbf{A}}^{[X-D]}
 \simeq A_{U}\otimes_{A}S_{\mathbf{P}}^{[X-D]}
 \simeq H^{0}(\mathbf{P}_{U},
 \omega_{\mathbf{P}/\mathbf{A}}^{n}(X))
 \simeq H^{0}(\mathbf{P}_{U},
 \omega_{\mathbf{P}/\mathbf{A}}^{n}(\log{X})).
 $$
 We assume that $0\leq p\leq n-2$.
 By Lemma~$\ref{c1}$ and Lemma~$\ref{c2}$, we have a commutative diagram
 $$
 \begin{CD}
 0&&0\\
 \downarrow&&\downarrow\\
 H^{0}(\mathbf{P}_{U},
 \omega_{\mathbf{P}/\mathbf{A}}^{n}((n-p-1)X))
 \otimes\Cl{(\mathbf{P}\smallsetminus E)}^{*}
 &\overset{\alpha_{1}}{\rightarrow}&
 H^{0}(\mathbf{P}_{U},
 \omega_{\mathbf{P}/\mathbf{A}}^{n}((n-p-1)X))\\
 \downarrow&&\downarrow\\
 \displaystyle{\bigoplus_{\rho\in\Sigma_{\pi}(1)}}
 H^{0}(\mathbf{P}_{U},
 \omega_{\mathbf{P}/\mathbf{A}}^{n}
 ((n-p-1)X+\mathbf{P}_{\rho}))
 &\overset{\alpha_{2}}{\rightarrow}&
 H^{0}(\mathbf{P}_{U},
 \omega_{\mathbf{P}/\mathbf{A}}^{n}((n-p)X))\\
 \downarrow&&\downarrow\\
 H^{0}(\mathbf{P}_{U},
 \omega_{\mathbf{P}/\mathbf{A}}^{n-1}((n-p-1)X))
 &\overset{\alpha_{3}}{\rightarrow}&
 H^{0}(\mathbf{P}_{U},
 \omega_{\mathbf{P}/\mathbf{A}}^{n}((n-p)X)\vert_{X})\\
 \downarrow&&\downarrow&&&&\\
 0&&0,
 \end{CD}
 $$
 where the exactness of the left vertical sequence is shown by
 Theorem~$\ref{eu}$.
 By Lemma~$\ref{resol}$ and the vanishing theorem (Corollary~$\ref{v}$),
 the cohomology group
 $H^{n-p-1}(\mathbf{P}_{U},
 \omega_{\mathbf{P}/\mathbf{A}}^{p+1}(\log{X}))$
 appear in the cokernel of the map $\alpha_{3}$.
 Since the cokernel of the map $\alpha_{2}$ is
 $A_{U}\otimes_{A}R_{X/\mathbf{A}}^{[(n-p)X-D]}$, we have to show that
 the map $\alpha_{1}$ is surjective.
 It follows from the assumption that the class
 $[X]\in\Cl{(\mathbf{P}\smallsetminus E)}$ is not divisible by
 the characteristic of $k$.
\end{proof}
Using Theorem~$\ref{main}$, we have a description for the cohomology of
relative logarithmic forms on $X$.
\begin{corollary}
 \renewcommand{\labelenumi}{$(\arabic{enumi})$\hspace*{5pt}}
 \begin{enumerate}
  \item For $0\leq p\leq\frac{n-3}{2}$,
	$$
	A_{U}\otimes_{A}R_{X/\mathbf{A}}^{[(n-p)X-D]}
	\simeq
	H^{n-p-1}(X_{U},\omega_{X/\mathbf{A}}^{p}),
	$$
	and for $p=\frac{n-2}{2}$, there is an exact sequence
	$$
	H^{\frac{n}{2}}
	(\mathbf{P}_{U},\omega_{\mathbf{P}/\mathbf{A}}^{\frac{n}{2}})
	\longrightarrow
	A_{U}\otimes_{A}R_{X/\mathbf{A}}^{[(\frac{n}{2}+1)X-D]}
	\longrightarrow
	H^{\frac{n}{2}}(X_{U},\omega_{X/\mathbf{A}}^{\frac{n}{2}-1})
	\longrightarrow0.
	$$
  \item If $\mathbf{A}$ is nonsingular, and
	$\mathcal{O}_{\mathbf{P}}(-E)$ is generated by global
	sections, then for $\frac{n}{2}\leq p\leq n-1$,
	$$
	A_{U}\otimes_{A}R_{X/\mathbf{A}}^{[(n-p)X-D]}
	\simeq
	H^{n-p-1}(X_{U},\omega_{X/\mathbf{A}}^{p}),
	$$
	and for $p=\frac{n-1}{2}$, there is an exact sequence
	$$
	0\longrightarrow
	A_{U}\otimes_{A}R_{X/\mathbf{A}}^{[(\frac{n+1}{2})X-D]}
	\longrightarrow
	H^{\frac{n-1}{2}}(X_{U},\omega_{X/\mathbf{A}}^{\frac{n-1}{2}})
	\longrightarrow
	H^{\frac{n+1}{2}}
	(\mathbf{P}_{U},\omega_{\mathbf{P}/\mathbf{A}}^{\frac{n+1}{2}}).
	$$
 \end{enumerate}
\end{corollary}
\begin{proof}
 This is proved by the exact sequence
 $$
 0\longrightarrow
 \omega_{\mathbf{P}/\mathbf{A}}^{p+1}
 \longrightarrow
 \omega_{\mathbf{P}/\mathbf{A}}^{p+1}(\log{X})
 \longrightarrow
 \omega_{X/\mathbf{A}}^{p}
 \longrightarrow0
 $$
 and the vanishing theorem (Corollary~$\ref{v}$ and the next
 proposition).
\end{proof}
\begin{proposition}
 Let $\mathbf{A}=\mathbf{A}^{m}$ be an affine space, let $\mathbf{P}$ be
 a nonsingular toric variety, let
 $\pi:\mathbf{P}\rightarrow\mathbf{A}$ be a logarithmically smooth
 proper equivariant morphism, and let $\mathcal{L}$ be an invertible
 sheaf on $\mathbf{P}$.
 If
 $\mathcal{H}om_{\mathcal{O}_{\mathbf{P}}}
 (\mathcal{L},\mathcal{O}_{\mathbf{P}}(-E))$
 is generated by global sections, then for $0\leq q\leq p-1$,
 $$
 H^{q}(\mathbf{P},
 \omega_{\mathbf{P}/\mathbf{A}}^{p}
 \otimes_{\mathcal{O}_{\mathbf{P}}}\mathcal{L})=0.
 $$
\end{proposition}
\begin{proof}
 We define a locally free $\mathcal{O}_{\mathbf{P}}$-module by
 $\mathcal{F}=\mathcal{H}om_{\mathcal{O}_{\mathbf{P}}}
 (\mathcal{L},\omega_{\mathbf{P}/\mathbf{A}}^{n-p}(-E))$.
 By the duality theorem \cite{h} for the morphism $\pi$, there is an
 isomorphism
 $$
 \Ext_{\mathcal{O}_{\mathbf{P}}}^{n+i}
 {(\mathcal{F},\pi^{!}\varOmega_{\mathbf{A}}^{m})}
 \simeq
 \Ext_{\mathcal{O}_{\mathbf{A}}}^{i}
 {(\mathbf{R}\pi_{*}\mathcal{F},\varOmega_{\mathbf{A}}^{m})}.
 $$
 Since
 $\pi^{!}\varOmega_{\mathbf{A}}^{m}\simeq
 \mathcal{O}_{\mathbf{P}}(-D-E)$ and
 $\varOmega_{\mathbf{A}}^{m}\simeq\mathcal{O}_{\mathbf{A}}$,
 we have
 $$
 \Ext_{\mathcal{O}_{\mathbf{P}}}^{n+i}
 {(\mathcal{F},\omega_{\mathbf{P}/\mathbf{A}}^{n}(-E))}
 \simeq
 \Ext_{\mathcal{O}_{\mathbf{A}}}^{i}
 {(\mathbf{R}\pi_{*}\mathcal{F},\mathcal{O}_{\mathbf{A}})}.
 $$
 There is a spectral sequence
 $$
 E_{2}^{i,j}=
 \Ext_{\mathcal{O}_{\mathbf{A}}}^{i}
 {(R^{-j}\pi_{*}\mathcal{F},\mathcal{O}_{\mathbf{A}})}
 \Longrightarrow
 \Ext_{\mathcal{O}_{\mathbf{A}}}^{i+j}
 {(\mathbf{R}\pi_{*}\mathcal{F},\mathcal{O}_{\mathbf{A}})}.
 $$
 By Corollary~$\ref{v}$, we have $E_{2}^{i,j}=0$ for $-j\geq n-p+1$.
 Hence
 $$
 H^{q}(\mathbf{P},
 \omega_{\mathbf{P}/\mathbf{A}}^{p}\otimes\mathcal{L})
 \simeq
 \Ext_{\mathcal{O}_{\mathbf{P}}}^{q}
 {(\mathcal{F},\omega_{\mathbf{P}/\mathbf{A}}^{n}(-E))}
 \simeq
 \Ext_{\mathcal{O}_{\mathbf{A}}}^{q-n}
 {(\mathbf{R}\pi_{*}\mathcal{F},\mathcal{O}_{\mathbf{A}})}=0
 $$
 for $n-q\geq n-p+1$.
\end{proof}

\bigskip
{\scshape Graduate School of Science, Osaka University,\\
Toyonaka, Osaka, 560-0043, Japan}\\
{\itshape E-mail address}:
{\ttfamily atsushi@math.sci.osaka-u.ac.jp}
\end{document}